\documentclass[12pt,twoside]{paper}
\usepackage{amssymb,amsmath}
\usepackage[all,cmtip]{xy}
\usepackage{epic}
\usepackage{url}
\usepackage[english]{babel}
\usepackage{lscape}

\usepackage{graphics}
\usepackage[pdftex]{graphicx}

\def\bee{\begin{equation}}
\def\eee{\end{equation}}
\def\Li{{\rm li}}

\def \BB{{\cal B}}
\begin{document}

\thispagestyle{empty}
\bigskip\bigskip
\centerline{    }
\centerline{\Large\bf Search for  primes of the form  $m^2+1$}
\bigskip\bigskip\bigskip
\centerline{\large\sl Marek Wolf}
\bigskip
\centerline{ \sf Group of Mathematical Methods in Physics, University of Wroc{\l}aw}
\centerline{\sf Pl.Maxa Borna 9, PL-50-204 Wroc{\l}aw, Poland}
\centerline{\small e-mail:  mwolf@ift.uni.wroc.pl}

\begin{center}
{\bf Abstract}\\
\end{center}

\begin{minipage}{12.8cm}
The results of the  computer hunt for the primes of the form
$q = m^2+1$  up to $10^{20}$ are reported.
The number of sign changes of the difference
$\pi_q(x) - \frac{C_q}{2}\int_2^x{du \over \sqrt{u}\log(u)}$ and the error term
for this difference
is investigated. The analogs of the Brun's constant and the Skewes number are
calculated. An analog of the B conjecture of Hardy--Littlewood is formulated.
It is argued that there is no Chebyshev bias for primes of the form $q=m^2+1$.
All encountered integrals we were able to express by the logarithmic
integral.
\end{minipage}

\bigskip\bigskip\bigskip


\bibliographystyle{abbrv}

\section{Introduction}

This paper is devoted to investigation of the set of prime numbers
\bee
Q=\{2,5,17,37,101,197, 257, 401, 577 \ldots\}
\label{set_Q}
\eee
given by the quadratic polynomial  $m^2+1$ and let $q_n$ denote the $n$-th prime of this form.
By  the conjecture E of Hardy and Littlewood \cite{Hardy_and_Littlewood}
the number\footnote{we adopt here the
convention that all  functions describing usual prime numbers will have the subscript $q$
when applied to primes of the form $q = m^2+1$} $\pi_q(x)$  of primes  $ q< x$ of the form
$q = m^2 + 1$ is  given by
\bee
\pi_q(x) \sim C_q  \frac {\sqrt{x}}{\log(x)},
\label{conj_E}
\eee
where
\bee
C_q=\prod_{p\geq 3}\biggl(1-\frac{(-1)^{(p-1)/2}}{p-1}\biggr) =
1.372813462818246009112192696727\ldots
\eee
The primes $q_n$ were investigated in the past  both theoretically and numerically.
One of the strongest theoretical results
is the theorem of H. Iwaniec
\cite{Iwaniec}, who proved that there exist infinitely many integers $m^2 + 1$ which are
2-almost-primes. In 1998 H. Iwaniec and J. Friedlander   \cite{Iwaniec-Friedlander} have
proved that there is an infinity
of primes given by the polynomial of two variables $m, n$  of the form $m^2+n^4$,
thus the case  of the polynomial of one variable $m^2+1$ is not covered by their theorem.
Quite recently there appeared  two papers by S. Baier and L. Zhao
\cite{Baier-2007, Baier-Zhao}, treating the more general problem of the
primes of the form $m^2+k$ in average over square-free parameter $k$ in appropriate
intervals.
In the computational part we should cite the  papers by Shanks \cite{Shanks, Shanks1959}
and Wunderlich \cite{Wunderlich}; in the last paper the table of $\pi_q(x)$
for $x<1.96\times 10^{14}$ is given.

By analogy with the case of all primes, where substitution  of the
logarithmic integral ${\rm Li}(x) = \int_2^x du/\log(u)~$   instead of $x/\log(x)$ gives
better approximation for $\pi(x)$,
we  recast the original Hardy and Littlewood  conjecture  E   in the form:
\bee
\pi_q(x) \sim {1 \over 2} C_q  \int_2^x {du \over \sqrt{u}\log{u}} 
\label{conj_calka}
\eee
\bee
= C_q\left({\sqrt{x} \over \log{x}} +
{2\sqrt{x} \over \log^2{x}} + {8\sqrt{x} \over \log^3{x}} + {48\sqrt{x}\over \log^4{x}}+ \ldots+
{2^{n-1} (n-1)!\sqrt{x} \over \log^{n}{x}}+\ldots\right).
\nonumber
\eee

The series in parenthesis above is asymptotic one: the terms initially decrease, but
for sufficiently large $n$ they become to increase. Thus this series has to be cut
at such $n_0$, which depends on $x$, that the  $n_0$-th term  (and consecutive terms)
is larger than previous one with $n_0(x)-1$ what gives for the threshold at which
the series (\ref{conj_calka}) should be cut the inequality:
\bee
n_0(x)> \frac{1}{2} \log(x)+1.
\label{cut}
\eee
Besides this asymptotic representation the integral (\ref{conj_calka}) can
be linked to the logarithmic integral by the change of variables $u=t^2$:  
\bee
\int_a^b {du \over \sqrt{u}\log{u}} = \int_{\sqrt{a}}^{\sqrt{b}} {dt \over \log{t}}
=\Li (\sqrt{b})- \Li (\sqrt{a}).
\label{calka_Li}
\eee
Here we use the following convention for the lower limit of integration:
\bee
\Li(x)=v.p. \int_0^x \frac{du}{\log(u)} = \lim_{\epsilon \rightarrow 0}\left(
\int_0^{1-\epsilon} \frac{du}{\log(u)} + \int_{1+\epsilon} ^x \frac{du}{\log(u)}\right)
\eee
Integration by parts gives the asymptotic expansion:
\bee
\Li(x) \sim 
\frac{x}{\log( x)}  + \frac{x}{\log^2(x)} + \frac{2x}{\log^3( x)} +
\frac{6x}{\log^4(x)} + \cdots \frac{n!x}{\log^{n+1}(x)}.
\label{Li_asymptotic}
\eee
which should be cut at $n_0=\log(x)$. There is a series giving $\Li(x)$ for
all $x$ and  quickly convergent which has $n!$ in denominator and $\log^n(x)$
in nominator instead of opposite order in (\ref{Li_asymptotic})
(see \cite[Sect. 5.1]{abramowitz+stegun})
\bee
\Li (x) =  \gamma +
\log \log(x) + \sum_{n=1}^{\infty} {\log^{n}(x)\over n \cdot n!}
\quad {\rm for} ~ x > 1 ~ ,
\label{Li-series}
\eee
Here $\gamma=0.57721566490153286...$ is the Euler-Mascheroni constant.
Even faster converging series was discovered by Ramanujan \cite[p.123]{BerndtIV}:
\bee
{\rm li} (x) = \gamma + \log( \log( x)) + \sqrt{x} \sum_{n=1}^{\infty}
\frac{ (-1)^{n-1} (\log( x))^n} {n! \, 2^{n-1}}
\sum_{k=0}^{\lfloor (n-1)/2 \rfloor} \frac{1}{2k+1} \quad {\rm for} ~ x > 1 ~ .
\label{Li_R}
\eee

Skipping the number $2=1^2+1$ all odd primes $q_n$ can be expressed by the form
$m^2+1$ only for $m$ even, thus we
put $m=2k$,  i.e. we are looking for primes of the form $4k^2+1$. We  have collected
in one  file of the size roughly 3.5 GB
values of $k$ (thus  the prime  $2=1^2+1$ is absent here!) up to $k=4999999978$ what corresponds to all primes of the form $m^2+1<10^{20}$.
The compressed data occupies 760 MB and is available for downloading from
 \url{http://www.ift.uni.wroc.pl/~mwolf/4k-2+1-data.zip}.
Initially up to $6.65\times 10^{16}$ I  have used my own program written in Fortran
for Alpha DEC workstation,
but to reach $10^{20}$ we  have used the free package PARI/GP \cite{Pari} developed especially
for number theoretical purposes. To scan the interval $(6.65\times 10^{16}, 10^{20})$
it took about 10 days of CPU time on the PC with the clock 2.5 GHz using
the  built-in PARI  very fast function $isprime(p)$. The Table I
gives comparison of $\pi_q(x)$ with conjectures (\ref{conj_E}) and (\ref{conj_calka}).
Among those 312,357,934 values of $k$ there were 11,864,645 such $k$ that they
in turn were primes,
i.e. there were 11,864,645 pairs of numbers $(p, 4p^2+1)$  both  being prime.
We have checked separately that  up to $10^{22}$  there are 96,817,209 such pairs
$(p, 4p^2+1)$ \cite{Wolf}.

As it is seen from the fourth
column, the actual number $\pi_q(x)$ is always larger than prediction (\ref{conj_E}).
Contrary to this the ratio of $\pi_q(x)$ to  the  integral (\ref{conj_calka}) sometimes
is larger than 1 and sometimes is smaller than 1, see  the last column. It means,
that the
difference $\pi_q(x)-\frac{1}{2}  C_q  \int_2^x {du \over \sqrt{u}\log{u}}$ changes
the sign,
see the Sect.3. In the next Section 2 we discuss the problem of the error term, in Sect.4
the
analog of the Brun's constant is calculated. Sections 5  and 6 contains some heuristics
about the
distribution of $k$ in $4k^2+1$ giving the prime. We formulate an analog of the
B conjecture of Hardy and Littleewood for the case of primes $q=m^2+1$.
In Sect. 6 we discuss the distribution of the gaps between consecutive $k's$ giving the
prime $4k^2+1$. Heuristic arguments allow us to make conjecture about the growth
of the difference $q_{n+1}-q_n=\mathcal{O}(\sqrt{q_n}\log^2(q_n))$.
The last Section  7 contains the discussion of the analog of the Chebyshev's Bias
for the case of primes in  $Q$.

\section{The problem of error term}

Nothing is known about the error term for the formula (\ref{conj_calka})
(see however \cite{ Baier-2007, Baier-Zhao}),
thus the only way to gain some information and intuition is to appeal to the available
computer data. The Figure 1 presents the plot of the difference
\bee
|\Delta_q(x)|  = \left| \pi_q(x) - {1 \over 2} C_q  \int_2^x {du \over \sqrt{u}\log{u}}\right|.
\eee
for $x \in (10, 10^{20})$.  As it is seen from this plot the graph of the error term
is very erratic, thus the plot of
the  maximal value of the  absolute difference:
\bee
\omega(x)=\max_{2<t<x}|\Delta_q(t)|,
\eee
what is a kind of envelope for $|\Delta_q(x)|$, is also  plotted  in green in the Figure 1.

The error term present in the ordinary Prime Number Theorem (PNT)
under the Riemann Hypothesis is $\sqrt{x}\log(x)$:
\bee
\pi(x) = {\rm Li} (x) + {\mathcal{O}}(\sqrt{x}\log(x)).
\eee
which can be written in the slightly weaker form
\bee
\pi(x) = {\rm Li} (x) + {\mathcal{O}}(x^{\frac{1}{2}+\epsilon}).
\eee

\vskip 0.4cm
\begin{center}
{\sf TABLE {\bf I}}\\
\bigskip
\begin{tabular}{||c|c|c|c||c|c||} \hline
$x$ & $ \pi_q(x) $ & $ C_q\sqrt{x}/\log(x)$ & $\pi_q(x)/$ eq. (\ref{conj_E}) & formula (\ref{conj_calka}) & $ \pi_q(x)$/eq. (\ref{conj_calka})\\ \hline
$10^{6}$ &        112 &              99 &       1.12713 &             122 &       0.91869 \\ \hline
$10^{7}$ &        316 &             269 &       1.17325 &             318 &       0.99440 \\ \hline
$10^{8}$ &        841 &             745 &       1.12847 &             855 &       0.98321 \\ \hline
$10^{9}$ &       2378 &            2095 &       1.13516 &            2357 &       1.00888 \\ \hline
$10^{10}$ &       6656 &            5962 &       1.11639 &            6610 &       1.00696 \\ \hline
$10^{11}$ &      18822 &           17140 &       1.09815 &           18787 &       1.00184 \\ \hline
$10^{12}$ &      54110 &           49684 &       1.08909 &           53971 &       1.00258 \\ \hline
$10^{13}$ &     156081 &          145028 &       1.07621 &          156386 &       0.99805 \\ \hline
$10^{14}$ &     456362 &          425861 &       1.07162 &          456405 &       0.99991 \\ \hline
$10^{15}$ &    1339875 &         1256912 &       1.06601 &         1340089 &       0.99984 \\ \hline
$10^{16}$ &    3954181 &         3726285 &       1.06116 &         3955222 &       0.99974 \\ \hline
$10^{17}$ &   11726896 &        11090399 &       1.05739 &        11726340 &       1.00005 \\ \hline
$10^{18}$ &   34900213 &        33122538 &       1.05367 &        34903278 &       0.99991 \\ \hline
$10^{19}$ &  104248948 &        99229889 &       1.05058 &       104251624 &       0.99997 \\ \hline
$10^{20}$ &  312357934 &       298102838 &       1.04782 &       312353427 &       1.00001 \\ \hline
\end{tabular} \\
\end{center}
\vskip 0.4cm

The $\sqrt{x}$ behavior is confirmed by the computer
data, see e.g. \cite[Table 14, p. 175]{Ribenboim} or \cite[Table 5 and 6]{Granville-Races},
where  the difference ${\rm Li}(x) - \pi(x)$  has roughly half digits of the value $\pi(x)$.
Because there are roughly $\sqrt{x}$ candidates for primes of the form $m^2+1$ up to
$x$ it is natural to expect that the error term for (\ref{conj_calka}) will be
square root of the error term for PNT. This heuristic seems to be confirmed by the
fact, that $\omega(x)$  is well approximated by the power--like error term:
\bee
\alpha_1 x^{\beta_1},~~~~~~~~\alpha_1= 0.38\ldots, ~~\beta_1=0.22\ldots .
\label{error}
\eee
and indeed here $\beta_1 \approx 1/4$.
This function was obtained by fitting the straight line to the points $\log(\omega(x))$
vs $\log(x)$ for $x>10^9$ by the least-square method and  to bound the difference
$\omega(x)$ from above it is sufficient to shift the above curve (\ref{error}) parallel
up to leave the  plot of $|\Delta_q(x)|$  below. In the Figure 1 the function
$5x^\beta_1$ was chosen, but
at least for $x<10^{20}$ the smaller choice for the constant hidden in the
big-$\mathcal{O}$ in $\omega(x)=\mathcal{O}(x^{\beta_1})$
will also do. These heuristic arguments and computer data lead us to guess
the following\\

{\bf Conjecture 1:}
\bee
\pi_q(x) = {1 \over 2} C_q  \int_2^x {du \over \sqrt{u}\log{u}} + \mathcal{O}(x^{\frac{1}{4}+\epsilon}).
\label{Conj_1}
\eee
More stringent error term $\mathcal{O}(x^{\frac{1}{4}}\sqrt{\log(x)})$ is also
a possibility which cannot be ruled out by available data. Let us notice, that
(\ref{Conj_1}) is heuristically supported by the relation (\ref{calka_Li}) and
$\pi(\sqrt{x}) =
\Li(\sqrt{x})+ \mathcal{O}(x^{1/4})$. Another support in favor of (\ref{Conj_1})
will be given in Sect.7 (see Fig.8).

The difference $\Delta_q(x)$ fluctuates roughly symmetrically  around zero.
As the computer check of possible future oscillations theorems we present in the
Fig.2 the plot of $\Delta_q(x)/x^{1/4}$. The amplitude of this function practically
does not change in the interval $(10^2, 10^{20})$ giving support in favor of the
error term $\mathcal{O}(x^{1/4})$.

\section{ An analog of the Skewes number for primes of the form $m^2+1$}

It turns out that the there is a lot of sign changes of the difference
\bee
\Delta_q(x)=\pi_q(x) - {1\over 2} C_q\int_2^x{du \over \sqrt{u}\log(u)}
\eee
in the investigated
interval $x \in (2, 10^{20}$). In the generic problem of all prime numbers it was shown
by J.E. Littlewood  in the 1914 \cite{Littlewood} (see also \cite{Ellison}) that
the difference between the number of primes smaller than $x$ and the
logarithmic integral li$(x)$ infinitely often changes the sign.
The smallest value
$x_S$ such that for the first time the difference $\Delta(x)=\pi(x) - \Li(x)$
changes the sign is called Skewes number
and the lowest present day known estimate of the Skewes number is around
$10^{316}$, see \cite{Bays} and \cite{Demichel}.
However in the case of the primes given by the quadratic polynomial $m^2+1$ the first
sign change of the difference $\Delta_q(x)$ occurs already at the prime
$q_{13}=2917=54^2+1$ and
there are 20634 such sign changes up to $10^{20}$.
Let $\nu_q(T)$ denotes the number of sign changes of the function $\Delta_q(x)$ for
$x\in(2,T)$. The Fig.3  presents the plot of the function $\nu_q(T)$ for
$T<10^{20}$.  The fitting of the power--like dependence of $\nu_q(T)$ on $T$ gives
parameters  which depend on the number  of discarded initial points. For example
fitting $\log(\nu_q(T))$ vs $\log(T)$ for $T\in (10^7, 10^{20})$ gives $\nu_q(T)\sim T^{0.23308}$
while for $T\in (10^{12}, 10^{20})$  we obtained $\nu_q(T)\sim T^{0.23834}$, and we have
checked  for other intervals of $T$ that the first digits $0.23$ persists, thus we
write:
\bee
\nu_q(T)\approx \alpha_2 T^{\beta_2}~~~~\alpha_2=0.38\ldots, ~~\beta_2=0.23\ldots
\label{Skewes}
\eee
Let us mention that for the case of  all primes
Knapowski \cite{Knapowski} proved that the number of sign changes of $\Delta(x)=$
in the interval $(1,T)$
\bee
\nu(T)\geq e^{-35}\log \log \log \log T
\label{nu}
\eee
provided $T\geq \exp\exp\exp\exp(35)$. There is a remarkable coincidence in the values of the
parameters $\alpha$ and $\beta$ present in the fits (\ref{error}) and (\ref{Skewes})
and when the functions  $\omega(x)$ and $\nu_q(T)$   are plotted on
the same graph they appear very close to each other,
despite the fact that $\omega(x)$ and $\nu_q(T)$ represent quantities at first sight
unrelated.

\section{ The analog of the Brun constant}

Because $\sum_{n=1}^\infty \frac{1}{n^2}=\pi^2/6<\infty$ thus the sum of reciprocals of
all primes of the form $q = m^2+1$ is trivially convergent:
\bee
\sum_{q\in Q} \frac{1}{q} =\frac{1}{2}+\frac{1}{5}+ \frac{1}{17}+ \frac{1}{37}+\frac{1}{101}+\ldots < \infty
\label{suma_1_nad_p}
\eee
but the actual numerical value of this sum is unknown \cite{Plouffe-email}.
 In  1919 Brun \cite{Brun} has shown that the sum of
reciprocals of all twin primes is finite:
\bee
\BB_2=\left( {1 \over 3} + {1\over 5}\right) + \left( {1 \over 5} + {1\over
7}\right) + \left( {1 \over 11} + {1\over 13}\right) + \ldots < \infty.
\label{def_B2}
\eee
Numerically $\BB_2=1.9021605823\ldots$, see \cite{Nicely_Brun}, \cite{Brent}. It is natural to
call the above sum (\ref{suma_1_nad_p}) the Brun's constant for primes of the form
$q = m^2+1$ and
denote it by $\BB_q$. From the computer data we  can calculate the finite size
approximations:
\bee
\BB_q(x) = \sum_{x>q\in Q} \frac{1}{q}.
\label{Brun_x}
\eee
From  the  integral test for convergence of the series in the  form:
\bee
\sum_{n=N}^\infty  f(n)\le f(N)+\int_N^\infty  f(u)\,du.
\eee
we  have:
\bee
\BB_q  =  \sum_{x>q\in Q} \frac{1}{q} + \sum_{x<q\in Q} \frac{1}{q}
<\BB_q(x) + \sum_{n^2>x} \frac{1}{n^2}
<\BB_q(x) + \frac{1}{x}+ \int_{\sqrt{x}}^\infty \frac{du}{u^2}=
\BB_q(x)+\frac{1}{x}+\frac{1}{\sqrt{x}}
\eee
Thus from this  trivial inequality we can expect for $x=10^{20}$ the accuracy of 10 digits
for $\BB_q$. Indeed, from the second  column of the Table II we see that the number
of stabilizing digits of $\BB_q(x)$ is roughly half of the digits of the exponent
of $x$ in the first column. However, using some heuristics it is possible to obtain
from $\BB_q(10^{20})$
15 digits of $\BB_q$.  Namely, we can obtain the analytical
formula for dependence of $\BB_q(x)$ on $x$. From the equation (\ref{conj_calka})
it follows that the chance to find
a prime of the form $m^2+1$ around $x$  is  $\frac{C_q}{2\sqrt{x}\log(x)}$,\footnote{Let us
remark, that in the case of all primes $\pi(x)\sim \int_2^x \frac{du}{\log(u)}=\frac{x}{\log(x)}
+\ldots$ as well as of twin primes $\pi_2(x)\sim C_2 \int_2^x \frac{du}{\log^2(u)}=
C_2\frac{x}{\log^2(x)}+\ldots$ the chance to find a prime around $x$ following from
dividing the first terms on r.h.s. by $x$ coincide with integrands on l.h.s, what is not true
in the case  of (\ref{conj_calka}), where the factor $\frac{1}{2}$ makes the
difference.}   thus we can write:
\bee
\BB_q(x) = \sum_{x>q\in Q} \frac{1}{q} \approx \BB_q(\infty)-\frac{C_q}{2}\int_x^\infty \frac{du}{u^{3/2}\log(u)}
\label{Brun_x_2}
\eee
Integrating by parts gives:
\bee
\int_x^\infty \frac{du}{u^{3/2}\log(u)} = -\frac{2}{\sqrt{x}  \log{x}} +
\frac{4}{\sqrt{x}\log^2{x}} 
\ldots +(-1)^n\frac{2^n (n-1)!}{\sqrt{x} \log^{n}{x}}+\ldots
\label{calka_3_2}
\eee
This series is asymptotic one and the condition for the dropped terms is
$n>\frac{1}{2}\log(x)$ --- the same as threshold (\ref{cut}). But by the change of the
variable $u=1/\sqrt{t}$  it is possible to express the above integral by the
logarithmic integral:
\bee
\int_a^b \frac{du}{u^{3/2}\log(u)} = \Li\left(\frac{1}{\sqrt{b}}\right) -
\Li \left(\frac{1}{\sqrt{a}}\right).
\eee
This formula is useless for our purposes, because (\ref{Li-series}) or (\ref{Li_R})
is valid for  $x>1$.    The analog of
(\ref{Brun_x_2}) for usual Brun's constant is given by
\bee
\BB_2(\infty)=\BB_2(x)+\frac{2C_2}{\log(x)},
\label{Brun}
\eee
where $C_2=\prod_{p>2}(1-\frac{1}{(p-1)^2}) = 0.66016\ldots $ is the Twins constant.

The third column in Table II gives the sample of values
\bee
\BB_q^\star(x) = \BB_q(x)+\frac{C_q}{2}\int_x^\infty \frac{du}{u^{3/2}\log(u)}
\label{suma}
\eee
which are supposed to be constant and equal to $\BB_q(\infty)$.
As it seen from the Table II indeed with  increasing $x$ growing number of digits of the
sum (\ref{suma}) is stabilizing. To produce the data for this Table we calculated
in PARI the value of $C_q$  with over 30 digits accuracy using the formula (10)
from \cite{Shanks1959}. We calculated the finite approximations $\BB_q(x)$ in
DEC Fortran using the quadruple precision (REAL*16) with 33 decimal digits and
the Mathematica v.7  was used to calculate  integrals
(\ref{calka_3_2}) with over 30 digits of accuracy.

\vskip 0.4cm
\begin{center}
{\sf TABLE {\bf II}}\\
\bigskip
\begin{tabular}{||c|c|c||} \hline
$x$ & $ \BB_q(x) $ & $ \BB_q^\star(x) $ \\ \hline  
$10^{2}$ &         0.79575154653780163  &          0.81798411614326626  \\ \hline
$10^{3}$ &         0.81119372199008314  &          0.81625948818529443  \\ \hline
$10^{4}$ &         0.81335296609757082  &          0.81460891435244917  \\ \hline
$10^{5}$ &         0.81432340308206016  &          0.81465046053349854  \\ \hline
$10^{6}$ &         0.81450766696668914  &          0.81459563341678427  \\ \hline
$10^{7}$ &         0.81457232824055968  &          0.81459653781164721  \\ \hline
$10^{8}$ &         0.81458971836435488  &          0.81459649673657666  \\ \hline
 \vdots & \vdots & \vdots \\ \hline
$10^{13}$ &         0.81459655805079256  &          0.81459657169340710  \\ \hline
$10^{14}$ &         0.81459656768254395  &          0.81459657170479096  \\ \hline
$10^{15}$ &         0.81459657051185242  &          0.81459657170321129  \\ \hline
$10^{16}$ &         0.81459657134862317  &          0.81459657170292293  \\ \hline
$10^{17}$ &         0.81459657159724819  &          0.81459657170299049  \\ \hline
$10^{18}$ &         0.81459657167131276  &          0.81459657170297237  \\ \hline
$10^{19}$ &         0.81459657169347024  &          0.81459657170297623  \\ \hline
$10^{20}$ &         0.81459657170012661  &          0.81459657170298816  \\ \hline
\end{tabular} \\
\end{center}
\vskip 0.4cm

It is seen from the last column that starting with $x=10^{16}$ all
first 14 digits remain the same ---
the change appears at the 14-th place after the dot. In the paper
\cite{Nicely} Nicely  has performed complicated statistical analysis to get the
95 \% confidence interval for the value of $\BB_2$.
In our case it is possible to estimate the error
appearing in (\ref{suma}) by using the form of the function $\omega(x)$
given by Conjecture 1. Namely, the  ``density'' of  the error for the chance
$\frac{C_q}{2\sqrt{x}\log(x)}$ to find the
prime of the form $m^2+1$ around $x$ is less than  $\mathcal{O}(x^{-3/4})$, thus we have:
\bee
|\BB_q(\infty)-\BB_q^\star(x)| = \mathcal{O} (x^{-3/4}).
\eee
From this we see, that for $x=10^{20}$ the value of $\BB_q(\infty)$ lies in the interval
of approximate length $10^{-15}$ around  $\BB_q^\star(10^{20})$ and we can claim
that with 15 digits accuracy
\bee
\BB_q(\infty)\equiv\BB_q=0.81459657170299\ldots.
\eee

S. Plouffe has checked using his Symbolic Inverse Calculator
(http://pi.lacim.uqam.ca/eng/),
that this  constant  can not be expressed by other mathematical constants
\cite{Plouffe-email}, thus
the value of $\BB_q$ could be treated as a new mathematical constant. The comparison
of  numbers in the second and third column reveals that addition of  the term
$\frac{C_q}{2}\int_x^\infty \frac{du}{u^{3/2}\log(u)}$ causes that about 3-4
digits more than in the values of  $\BB_q(x)$ alone settle down, thus the rate
of convergence of $\BB_q^\star(x)$ is a few orders faster than that of $\BB_q(x)$.

The correctness of the choice of $\frac{1}{2}$ in front of the integral in
(\ref{Brun_x_2}) can be checked by comparing the values of the equation
\bee
\BB_q(x_2)-\BB_q(x_1)=\frac{C_q}{2}\int_{x_1}^{x_2} \frac{du}{u^{3/2}\log(u)}
\label{test}
\eee
following from (\ref{Brun_x_2}) with the actual computer data. For example for
$x_1=10^6$ and $x_2=10^{20}$ for the l.h.s. of (\ref{test}) from the computer
data we  get  $0.000177\ldots$ while the r.h.s. is equal to $0.000176\ldots$.

Let us mention that two first primes from $Q$  give contribution
$\frac{1}{2} + \frac{1}{5} = 0.7$  to $\BB_q$, i.e. 86\% of the total
value $ 0.8145966\ldots$!

We can define the  Brun's measure  of the set of numbers $S=\{a_1,a_2,\ldots,\}$ as
the sum
\bee
\mathcal{M_B(S)}=\sum_i \frac{1}{a_i}
\eee

\noindent provided it is finite. Thus we can say that the Brun's measure  of the set of
twin primes is $1.9021605823\ldots$, while the Brun's measure  of the set
$Q$ is $\mathcal{M_B}(Q)=0.8145965717\ldots.$

\section{Analog of the conjecture B of Hardy--Littlewood}

The conjecture B of Hardy--Littlewood \cite{Hardy_and_Littlewood}    says that the number
\bee
\pi(d;  x)= \sum_{\substack{p<x\\ p-p'=d}} 1 = \mid \{ p<x: p  {~~and ~~} p-d { ~~ are~ primes}\} \mid
\eee
of prime pairs $p-d, p <x$ separated by $d$ ($d=2, 4, 6, \ldots$) is given by
\bee
\pi(d; x)\sim C_2 \prod_{\substack{p >2 \\p\mid d }}\left(\frac{p-1}{p-2}\right)\int_2^x \frac{dt}{\log^2(t)},
\label{H-L-B}
\eee
where
\bee
C_2 = 2\prod_{p >2} \left(1- \frac{1}{(p-1)^2}\right) = 1.320323631693739\ldots
\eee
is the Twin constant.  The integral in (\ref{H-L-B}) again can be expressed by the
logarithmic integral and thus calculated quickly from the series (\ref{Li-series})
or (\ref{Li_R}):
\bee
\int_a^b \frac{dt}{\log^2(t)} = \int_a^b \frac{dt}{\log(t)} + \frac{a}{\log(a)}
- \frac{b}{\log(b)} = \Li(b)- \Li(a) + \frac{a}{\log(a)} - \frac{b}{\log(b)}
\eee

The values of the gaps between  primes  $m^2+1$   and  $(m-d)^2+1$ grow linearly
with $m$, but we can formulate an analog of the conjecture B of  Hardy--Littlewood
when  we will focus on the gaps between $k'$s appearing in $4k^2+1$. Thus let us
define:
\bee
\pi_q(d;  x)=  \mid \{ 2k<\sqrt{x}:  4k^2+1 { ~~and ~~} 4(k-d)^2+1 { ~~ are~ both~ primes}\} \mid
\label{pi_q_d_x}
\eee
In contrast to all primes here $d=1, 2, 3, \ldots$.
The Fig. 4 presents the plots of $\pi_q(d; x)$
obtained from our computer data for $x=10^{10}, x=10^{12}, \ldots x=10^{20}$ and
$d\leq 300$.

There is a heuristic procedure  of Bateman and Horn \cite{Bateman} (see also
\cite[chap. 3]{Riesel})
allowing to guess the formula for the number of constellations of primes of different
types.
Let ${\bf f}=\{f_1(x), f_2(x), \ldots, f_l(x)\}$  be the set of distinctive irreducible
polynomials with
integral  coefficients and positive leading coefficient such, that
$f(x)=f_1(x) f_2(x) \ldots f_l(x)$ has no fixed
divisor $>1$. Let $\pi( {\bf f}; x)$ denote the number of positive integers $n<x$ such,
that all $f_1(n), f_2(n), \ldots, f_l(n)$ are simultaneously primes. Then the
Bateman--Horn conjecture reads:
\bee
\pi({\bf f}; x) \sim \prod_p \frac{1-\frac{w(p)}{p}}{(1-\frac{1}{p})^l}\int_2^x \frac{du}{\prod_{i=1}^l\log(f_i(u))}, 
\eee
where $w(p)$ is the number of {\it distinct} solutions to $f_1(u)f_2(u)\cdots f_k(u)\equiv 0
\pmod p$ with $u\in\{0,1,\ldots, p-1\}$. 
In our case $f_1(u)=4u^2+1, f_2(u)=4(u-d)^2+1$.
It is well known  \cite[chap. VII]{H-W} that equation $4u^2+1 \equiv 0  \pmod p$  has two solutions
for primes of the form $p=4n+1$ and does not have solutions for primes of the form
$p=4n+3$. For the case $p=4n+1$ if  $p\mid d$  then equations
$f_1(u)\equiv 0  \pmod p$ and $f_2(u)\equiv 0  \pmod p$ have the
same solutions, thus $w(p)=2$. When $p\nmid d$  then there are two possibilities:
$w(p)=3$ when $f_1(u)\equiv 0  \pmod p$ and $f_2(u)\equiv 0  \pmod p$ have one common solution,  or $w(p)=4$ when $f_1(u)\equiv 0  \pmod p$
and $f_2(u)\equiv 0  \pmod p$ have distinct solutions. The equations
$ f_1(u) \equiv 0  \pmod p$ and $f_2(u)\equiv 0  \pmod p$  can have one common
solution only when $p\mid \pm 2u^\star+d$, where $u^\star$ is the solutions of
$$4u^2+1 \equiv 0  \pmod p.$$
The solutions of this equation are of the form (see \cite[p.88]{H-W})
$$2u_{1,2}=\pm u^\star,    ~~~~u^\star = \left(\frac{p-1}{2}\right)!$$
Because the conditions $p\mid \pm 2u^\star+d$  for the case $w(p)=3$  can be written as  the one condition $p\mid d^2 - 4u^\star$ and  there is an identity
$4{u^\star}^2+1 \equiv 0  \pmod p$  we can write the condition for $w(p)=3$ simply as $p\mid d^2 +1$.
Thus we  have for $p=4n+1$

\bee
w(p)=\begin{cases}
4 &  \text{if~~} p\nmid d,~~p\nmid d^2+1   \\
3 &  \text{if~~} p\mid d^2+1 \\
2 &  \text{if~~} p\mid d
    \end{cases}
\eee

For $p=2$ and for $p=4n+3$ there are no solutions, i.e.  $w(p)=0$, hence finally the
product appearing in the Bateman--Horn conjecture takes the form:
\begin{eqnarray}
\hskip-5cm  \prod_p \frac{1-w(p)/p}{(1-1/p)^2}  = ~~~~~~~~~~~~~~~~~~~~~~~~~~~~~~~& \\
4 \prod_{p \equiv 3 \!\!\!\! \pmod 4}\frac{1}{(1-1/p)^2}
\prod_{\substack{p \equiv 1\!\!\!\!\! \pmod  4\\p\nmid d,~ p \nmid d^2+1 }}\frac{1-4/p}{(1-1/p)^2}
\prod_{\substack{p \equiv 1 \!\!\!\! \pmod 4\\p\mid d }}\frac{1-2/p}{(1-1/p)^2}
\prod_{\substack{p \equiv 1 \!\!\!\! \pmod 4\\p\mid d^2+1 }}\frac{1-3/p}{(1-1/p)^2}.  &\nonumber
\end{eqnarray}
We can get rid of the product over $p\nmid d$ by extending it to  the product over all
$p \equiv 1 \pmod 4 $ and simultaneously by dividing by an appropriate term. Because
the conditions $p \mid d$ and  $p\mid  d^2+1$  cannot be satisfied simultaneously
these additional factors can be incorporated into
the last two products above. Finally the factor describing oscillations
takes the form (we separated 4 to cancel it later with 4 coming from the degrees of
$f_1(u)$  and $f_2(u)$):
\bee
\prod_p \frac{1-w(p)/p}{(1-1/p)^4}= 4 C_1 P(d),
\label{iloczyn}
\eee
where the constant $C_1$:
\bee
C_1=\prod_{p \equiv 3\!\!\!\! \pmod 4}\frac{p^2}{(p-1)^2}\prod_{p \equiv
1\!\!\!\! \pmod 4}\frac{p(p-4)}{(p-1)^2} = 0.975245556223143537223292783\ldots
\label{constant-C_1}
\eee
and  $P(d)$ denotes the product:
\bee
P(d)=\prod_{\substack{p \equiv 1 \!\!\!\! \pmod 4\\p\mid d }}\frac{p-2}{p-4}
\prod_{\substack{p \equiv 1 \!\!\!\! \pmod 4\\p\mid d^2+1 }}\frac{p-3}{p-4}.
\label{product_P}
\eee

In above expressions  the  condition $p\equiv 1 \pmod 4$ means that products are
over primes $p\geq 5$, thus all these products are positive. In (\ref{product_P})
the conditions   $p\mid d$ and   $p\mid d^2+1$ are fulfilled only by finite number
of $p'$s, hence it is obvious that these products  are convergent.

Finally we obtain the number of such $k<x$ that both $f_1(k) = (2k)^2+1$ and
$f_2(k) = (2(k-d))^2+1$ are prime:
\bee
\pi(f_1, f_2; x)= C_1 \prod_{\substack{p \equiv 1 \!\!\!\! \pmod 4\\p\mid d }}\frac{p-2}{p-4}
\prod_{\substack{p \equiv 1 \!\!\!\! \pmod 4\\p\mid d^2+1 }}\frac{p-3}{p-4}
 \int_1^{x} \frac{du}{\log^2(2u)}
\eee
We put here for
a while 1 as the lower limit of integration, since $k=1$ gives the prime
$4\cdot 1^2+1=5$ (let us notice, that Ramanujan often did not specify
the lower limit of integration, see \cite[p.123] {BerndtIV}). Because we skip 1
in $m^2+1$
in  manipulations of integrals below, alternatively we can say that the lower
limit of integration is $\sqrt{5/4}=1.118033989\ldots$.

Usually we are interested  directly in the number of primes $4k^2+1<x$ and after
the change of the integration variable $u=\sqrt{t}/2$ we have
\bee
\int_{\sqrt{a}/2}^{\sqrt{b}/2} \frac{du}{\log^2(2u)}=
\int_a^b \frac{dt}{\sqrt{t}\log^2(t)}
\eee
and finally for the quantity $\pi_q(d; x)$  defined in (\ref{pi_q_d_x}) we
obtain
\bigskip
{\bf Conjecture 2:}
\bee
\pi_q(d; x)=  C_1 \prod_{\substack{p \equiv 1 \!\!\!\! \pmod 4\\p\mid d }}\frac{p-2}{p-4}
\prod_{\substack{p \equiv 1 \!\!\!\! \pmod 4\\p\mid d^2+1 }}\frac{p-3}{p-4}
\int_5^x \frac{dt}{\sqrt{t}\log^2(t)}
+~~error~~term
\label{conjecture2}
\eee
We cheated a little here replacing 4 by 5 as the lower limit of integration;
better possibility as the lower limit of integration is perhaps
the Soldner-Ramanujan constant $\mu=1.45136923488338105\ldots$ defined
by $\Li(\mu)=0$, see \cite[p.123]{BerndtIV}.
Appearing here integral by the change of the variable $t=u^2$ can be
expressed by the logarithmic integral:
\bee
\int_a^b \frac{dt}{\sqrt{t}\log^2(t)} =\frac{1}{2}(\Li(\sqrt{b}) - \Li(\sqrt{a})) +
\frac{\sqrt{a}}{\log(a)} - \frac{\sqrt{b}}{\log(b)}.
\eee

The  product $P(d)$ (\ref{product_P}) is responsible for characteristic oscillations
seen in the Fig. 4.
The conjecture  2 agrees with the computer data quite well. Instead of producing
some table to corroborate this  statement, we give ``visual argument'': in the Fig.4
by the red  boxes are plotted values of the quotient $\pi_q(d; x)/P(d)$ for $x=10^{12},
\ldots x=10^{20}$  and $d\leq 300$. In green are plotted values of the integral
$\int_2^x du/\sqrt{t}\log^2(t)$  for the same values of $x$.

\section{Heuristics on the gaps between adjacent $k$}

It is interesting to restrict the analysis from the previous section to the case of
consecutive values of $k$ giving the prime $4k^2+1$.  Let us define the quantity:
\[
h(d; x)= \{number~ of~ pairs~~ k < k^\prime=k+d, ~ such~ that~~
4k^2+1~ and~ 4{k^{\prime 2}} +1 < x ~ are ~
\]
\bee
consecutive ~ primes  ~ of~ the~ form~ m^2+1\}  =  \sum_{\substack{q_n<x\\ q_n-q_{n-1}=4d(\sqrt{q_n-1}-d)}} 1 .
\eee

The Fig. 5 presents the plot of $h(d; x)$ 
obtained from our computer data for $x=10^{10}, x=10^{12}, \ldots x=10^{20}$. On the
semi-logarithmic scale the points display characteristic
oscillations around straight lines representing the  fits to exponential
decrease  obtained by the
least-square method. To describe the oscillations we take the product $P(d)$ from
(\ref{product_P}).  The Fig.5 presents
the corroboration of this mechanism of oscillations: dividing the values of $h(d,x)$ by
the product $P(d)$ leaves pure exponential  decrease --- small  deviations could be
attributed to the fluctuations of $h(d,x)$ and incorporated into the (unknown) error term.

The Figure 5 suggests the following\\

{\bf Ansatz:}
\bee
h(d,x)=C_1 P(d) B(x) e^{-d A(x)}
\label{guess}.
\eee
The functions $A(x)$ and $B(x)$, giving the slopes and the
intercepts of straight lines seen in the Fig. 5,
can be determined by exploiting two selfconsistency conditions that
$h(d, x)$ has to obey just from the definition. First of all, the number of
all gaps between $k's$ is by one
smaller than the number of primes of the form $4k^2+1$ smaller than $x$:
\bee
\sum_{d=1}^{K(x)} h(d,x)=\pi_q(x)-1,
\label{identity_1}
\eee
where $K(x)$ denotes the largest gap between two consecutive $k,k^\prime <x$.
The second selfconsistency condition comes from the observation, that the
sum of distances between adjacent $k$ is equal to the $k$ producing the
largest prime $q = 4k^2+1 \leq x$. For large $x$ we can write:
\bee
\sum_{d=1}^{K(x)}  h(d,x) d=\frac{\sqrt{x}}{2}.
\label{identity_2}
\eee
The erratic behavior of the product $P(d)$ is an obstacle in calculation of the above
sums (\ref{identity_1}) and  (\ref{identity_2}). Thus we  will replace $P(d)$ by the mean value:
\bee
\sum_{d=1}^n P(d)=n s +  E(n), 
\eee
where  we  assume that the  unknown error term $E(n)$ is an increasing function of $n$
which grows slower than $n$:
\bee
\lim_{n\rightarrow\infty} \frac{E(n)}{n}=0
\eee
what means that:
\bee
s=\lim_{n\rightarrow\infty}\frac{1}{n}\sum_{d=1}^n \left(\prod_{\substack{p \equiv 1 \!\!\!\! \pmod 4\\p\mid d }}\frac{p-2}{p-4}
\prod_{\substack{p \equiv 1 \!\!\!\! \pmod 4\\p\mid d^2+1 }}\frac{p-3}{p-4}\right).
\eee
There is known  at least one example of  the error term which grows slower than $n$
for the similar problem.
Namely, E. Bombieri and H. Davenport \cite{Bombieri} have proved that
the number  $1/\prod_{p > 2}( 1 - {1 \over (p - 1)^2})$ is the arithmetical average
for the product $\prod_{p\mid d} \frac{p-1}{p-2}$:
\bee
\sum_{d=1}^n \prod_{p\mid d,p>2}{p-1\over p-2} =
{n \over \prod_{p > 2}( 1 - {1 \over (p - 1)^2})} + \mathcal{O}(\log^2(n)).
\eee

To get rid of $P(d)$ in (\ref{identity_1}) and  (\ref{identity_2}) the Abel summation
formula can be used in the form:
$$
\sum_{i=1}^n a_i b_i = -\sum_{i=1}^{n-1} S(i)c_i + A(n) b_n,
$$
where  $S(i)=a_1+ \dots a_i$  and $c_i=b_{i+1}-b_i$.
Putting here $a_i=P(i)$, $b_i=f(i)$ and replacing $E(1) < E(2) < \ldots < E(n-1)$ by
larger $E(n)$  we obtain:
\bee
\sum_{l=1}^n P(l) f(l) = s \sum_{l=1}^{n-1} f(l) + \mathcal{O}(f(n)E(n)).
\label{rownosc}
\eee
In our case $f(l)=C_1B(x)e^{-A(x)l}$ for equation (\ref{identity_1}) and $f(l)=C_1lB(x)e^{-A(x)l}$
for equation (\ref{identity_2}).  From the Fig.5 we see, that $h(d,x)$ decreases
exponentially with $d$ and to solve (\ref{identity_1}) and (\ref{identity_2}) for
$e^{-A(x)}$ and $B(x)$ we have to  drop the term $\mathcal{O}(f(n)E(n))$.
The sums (\ref{identity_1}) and  (\ref{identity_2})  are the
geometrical  and differentiated geometrical series respectively; because $h(d,x)$
decreases exponentially with $d$ we  have replaced $K(x)$ in (\ref{identity_1}) and
(\ref{identity_2})  by $\infty$ and (\ref{identity_1}),(\ref{identity_2})
turn into the equations:
\bee
\frac{sC_1 B(x) e^{-A(x)}}{1-e^{-A(x)}}=\pi_q(x),
\eee
\bee
\frac{sC_1 B(x) e^{-A(x)}}{(1-e^{-A(x)})^2}=\frac{\sqrt{x}}{2}.
\eee
The solutions for $e^{-A(x)}$ and $B(x)$  of the above equations are:
\bee
e^{-A(x)}=1 - \frac{2\pi_q(x)}{\sqrt{x}},m 
\label{wzor_A}
\eee
\bee
B(x)=\frac{2\pi_q^2(x)}{sC_1\sqrt{x}(1-\frac{2\pi_q(x)}{\sqrt{x}})}. 
\label{wzor_B}
\eee
Hence we have finally:
\bee
h(d;x)=\frac{2P(d)\pi_q^2(x)}{s\sqrt{x}}\left(1-\frac{2\pi_q(x)}{\sqrt{x}}\right)^{d-1}
~~+~~error~term.
\label{h_d_x_1}
\eee

The Table III gives a comparison of the formulas (\ref{wzor_A}) and (\ref{wzor_B})
with the slopes and intercepts of $\log(h(d;x)/C_1P(d))$ vs $d$ obtained from the
computer  data
by means of the least square method (roughly 1/4 values of $h(d;x)$ for largest $d$
were discarded to avoid the fluctuations of $h(d;x)$ in the region of large $d$).


\vskip 0.4cm
\begin{center}
{\sf TABLE {\bf III}}\\
\bigskip
\begin{tabular}{||c|c|c|c||c|c|c||} \hline
$x$ & $ e^{-A(x)} $ & $ 1-\frac{2\pi_q(x)}{\sqrt{x}} $ & $  e^{-\frac{2\pi_q(x)}{\sqrt{x}}} $ &  $ B(x)  $  &  $ \frac{1}{sC_1}\frac{2\pi_q^2(x)}{\sqrt{x}-2\pi_q(x)} $ & $ \frac{2\pi_q^2(x)}{sC_1\sqrt{x}}  $  \\ \hline
$10^{15}$ &  $   0.91014186 $ & $     0.91525886 $ & $      0.91875009 $  & $         75421.77  $ &  $     65826.02  $  &  $       60247.85  $   \\ \hline
$10^{16}$ &  $   0.91611160 $ & $     0.92091638 $ & $      0.92396266 $  & $        206903.63  $ &  $    180179.71  $  &  $      165930.45  $   \\ \hline
$10^{17}$ &  $   0.92141730 $ & $     0.92583260 $ & $      0.92851624 $  & $        570171.11  $ &  $    498478.78  $  &  $      461507.91  $   \\ \hline
$10^{18}$ &  $   0.92632759 $ & $     0.93019957 $ & $      0.93257992 $  & $       1573509.35  $ &  $   1389610.28  $  &  $     1292614.89  $   \\ \hline
$10^{19}$ &  $   0.93060739 $ & $     0.93406718 $ & $      0.93619375 $  & $       4388166.12  $ &  $   3904615.30  $  &  $     3647172.98  $   \\ \hline
$10^{20}$ &  $   0.93444100 $ & $     0.93752841 $ & $      0.93943975 $  & $      12318241.09  $ &  $  11044184.29  $  &  $    10354236.57  $   \\ \hline
\end{tabular} \\
\end{center}
\vskip 0.4cm

The expression for the mean value $s$ of the product $P(d)$ can be
obtained in the following way: Because the pairs of primes of the
form $(2k)^2+1$ and $(2k+2)^2+1$  (Shanks calls in \cite{Shanks1960ngt}  such pairs
Gaussian Twins) correspond to  $d=1$ thus they are  necessarily consecutive
$(q_n, q_{n+1})$  and the formula for the number $h(1,x)$ of
such pairs smaller than $x$ obtained from (\ref{h_d_x_1}) has to be
equal to $\pi_q(1;x) $ from (\ref{conjecture2}):
\bee
h(1;x) = \frac{2\pi_q^2(x)}{s\sqrt{x}}=C_1 \int_5^x \frac{dt}{\sqrt{t}\log^2(t)}
\label{tozsamosc}
\eee
For  large $x$ we have $\pi_q(x)= C_q\sqrt{x}/\log(x)$  and
\[
\int_5^x \frac{dt}{\sqrt{t}\log^2(t)}=\frac{1}{2}\Li(\sqrt{x})-\frac{\sqrt{x}}{\log(x)}=
\frac{2\sqrt{x}}{\log^2(x)} + \ldots
\]
where we have used two first terms of the asymptotic expansion (\ref{Li_asymptotic})
and fortunately the term $\sqrt{x}/\log(x)$ cancels out leaving on both sides of
(\ref{tozsamosc}) the same dependence on $x$
and  thus we obtain $s=C_q^2/C_1$:
\bee  s=\lim_{n\rightarrow\infty}\frac{1}{n}\sum_{d=1}^n
\left(\prod_{\substack{p \equiv 1 \!\!\!\! \pmod 4\\p\mid d }}\frac{p-2}{p-4}
\prod_{\substack{p \equiv 1 \!\!\!\! \pmod 4\\p\mid d^2+1
}}\frac{p-3}{p-4}\right)~ = ~ \frac{C_q^2}{C_1}.
\eee
In the full (mysterious) form it reads:
\bee
\lim_{n\rightarrow\infty}\frac{1}{n}\sum_{d=1}^n
\left(\prod_{\substack{p \equiv 1 \!\!\!\! \pmod 4\\p\mid d }}\frac{p-2}{p-4}
\prod_{\substack{p \equiv 1 \!\!\!\! \pmod 4\\p\mid d^2+1
}}\frac{p-3}{p-4}\right) =
\frac{\left(\prod_{p\geq 3}\biggl(1-\frac{(-1)^{(p-1)/2}}{p-1}\biggr)\right)^2}
{\prod_{p \equiv 3\!\!\!\! \pmod 4}\frac{p^2}{(p-1)^2}\prod_{p \equiv
1\!\!\!\! \pmod 4}\frac{p(p-4)}{(p-1)^2}}
\label{srednia}
\eee

We  were not able to prove this identity analytically --- usual methods
of  calculating the sums of arithmetic  functions, see e.g.
\cite[Chap.1]{Iwaniec-Kowalski}, are not applicable here because
$P(d)$ is not a multiplicative function. The computer checking of
(\ref{srednia}) for large number of terms on the l.h.s. is also
difficult because the calculation of the average has almost cubic complexity in $n$
(i.e. finding the value of l.h.s. of  (\ref{srednia})  involves $\mathcal{O}(n^3/\log^2(n))$ operations).
We  have calculated the sum of $P(d)$ for $d$ up to 150000 and we obtained
$\sum_{d=1}^{150000} P(d) / 150000 = 1.93242674\ldots$, while the
r.h.s. of (\ref{srednia}) is   $(1.3728134628)^2/0.97524552 =1.93245368$,
thus the first 5 digits are the same.

In \cite{Shanks1960ngt}, \cite[p. 90]{Finch}, \cite{Sebah} heuristically the formula
for  $h(1;x)$ was obtained in the form:
\bee
h(1, x) = F\frac{\sqrt{x}}{\log^2(x)}~~+~~error~term
\eee
where:
\bee
F=\frac{\pi^2}{2}\prod_{p \equiv 1 \!\!\!\! \pmod
4}\left(1-\frac{4}{p}\right)\left(\frac{p+1}{p-1}\right)^2 =
1.9504911124462870744465855658\ldots.
\eee
In \cite{Sebah}  the 50 digits of this constant are given. Therefore we have $C_1=F/2, ~~
s=2C_q^2/F $
and the combination on the r.h.s of (\ref{srednia}) can be transformed to the form:
\bee
s=\frac{2C_q^2}{F}=\frac{1}{4}\prod_{p \equiv 1\!\!\!\! \pmod 4}\frac{(p-2)^2(p-1)^2}{p(p-4)(p+1)^2}
\prod_{p \equiv 3\!\!\!\! \pmod 4} \left(\frac{p+1}{p-1}\right)^2
\eee

Finally we state the\\

{\bf Conjecture 3:}
\bee h(d,x)= \frac{F P(d)\pi_q^2(x)} {C_q^2\sqrt{x}}
\left(1-\frac{2\pi_q(x)}{\sqrt{x}}\right)^{d-1}~~~+~~error~term.
\label{conjecture3}
\eee
For large $x$ we can simplify considerably the above formulas by  writing
$e^{-A(x)}=1-\frac{2\pi_q(x)}{\sqrt{x}}$ as
$e^{-A(x)}=e^{-\frac{2\pi_q(x)}{\sqrt{x}}}$ and $B(x)=\frac{2\pi_q^2(x)}{sC_1\sqrt{x}}$,
therefore in the limit of large $x$ we have:
\bee
h(d,x)= \frac{F P(d)} {C_q^2\sqrt{x}}
\pi_q^2(x) e^{-d \frac{2\pi_q(x)}{\sqrt{x}}}~~~+~~error~term.
\label{conjecture3b}
\eee
The Table III gives the comparison of the quantities $e^{-A(x)}$ and $B(x)$ obtained
from least-square method applied to $\log(h(d; x)/C_1 P(d))$ vs $d$ and analytical
expressions for them.

As a corroboration of the above conjectures
we will obtain the formula for the maximal gap $K (x)$ between two
consecutive values of $k<\sqrt{x}/2$ giving the prime $4k^2+1$.
Assuming, that the maximal gap
$K(x)$ appears only once we have the
equation $h(K(x),x)=1$ and putting the Hardy--Littlewood formula
for $\pi_q(x)$  in (\ref{conjecture3b}) and replacing $P(d)$ by $s$  we obtain:
\bee K(x) \sim
\frac{\log(x)}{2C_q}\left(\frac{1}{2}\log(x)+\log(2C_q^2)-2\log(\log(x))\right)
\eee
what for large $x$ goes to the
\bee
K(x) \sim \frac{1}{4C_q}\log^2(x).
\label{max_K}
\eee
Because for all primes it is widely believed that $G(x)\equiv{\max_{(p_n, p_{n-1}<x)}}
 (p_n-p_{n-1})
\sim \log^2(x)$, we see that $K(x)$ differs from $G(x)$ just by a constant (but here only
a fraction of $k's$ are primes!).
The largest gap between adjacent $k$ giving  $4k^2+1<10^{20}$ was 290. The comparison
of this formula with computer
data is shown in the Figure 6. From (\ref{max_K}) we deduce the following\\

{\bf Conjecture 4:}
\bee
q_{n+1}-q_n=\mathcal{O}(\sqrt{q_n}\log^2(q_n))
\label{Conjecture4}
\eee
Let us recall that for all primes the Riemann Hypothesis gives
$p_{n+1}-p_n=\mathcal{O}(p_n^{\frac{1}{2}+\epsilon})$ for any $\epsilon >0$, but
in reality gaps between consecutive primes are smaller and the Cramer conjecture
\cite{Cramer} states that $p_{n+1}-p_n=\mathcal{O}(\log^2(p_n))$, see
however  \cite{Granville}. Because our Conjecture 4 is obtained from
the guessed formula for maximal gap between $k's$ we expect (\ref{Conjecture4}) to be close
to the optimal bound for $q_{n+1}-q_n$.

\section{Analog of the Chebyshev's  bias}

For ordinary primes  the Dirichlet's Theorem  on the primes in
arithmetical progressions
asserts  that the number $\pi(x; 4, 1)$ of primes $<x$ giving 1
as the remainder when divided by 4 should be equal to  the number $\pi(x; 4, 3)$ of
primes  $<x$ giving 3 as the remainder when divided by 4, see e.g. \cite{Rubinstein_Sarnak}.
But the  direct inspection
shows that for small $x$ there are more primes $p\equiv 3 \pmod 4$ than $p\equiv 1 \pmod 4$:
$\pi(x; 4, 3)>\pi(x; 4, 1)$,
what is called  Chebyshev's  bias \cite{Kaczorowski}.
For the first time 1's takes the lead at $p_{2946}=26861$ --- up to this prime
3's win or there is a tie: $\pi(x; 4, 3)\geq \pi(x; 4, 1)$ for $x<26861$.
The next time  $\pi(x; 4, 3)< \pi(x; 4, 1)$ at 616841 and in general there
is preponderance of primes in the progression $4n+3$.  The same
phenomenon was observed in  other arithmetical
progressions \cite{Granville-Races}. However for the case of number of twins $\pi(2;x)$
(primes pairs ${p,p+2}$) and number of  cousins $\pi(4; x)$
(primes pairs ${p,p+4}$) which by the B conjecture of Hardy and Littlewood
(\ref{H-L-B})
 should be the same $\pi(2; x)\approx \pi(4; x)$,
there is no Chebyshev bias at least up to $2^{42}\approx 1.4\times 10^{12}$: sometimes
twins and
sometimes cousins take the lead,  \cite{Wolf1998}. In this paper it was shown
numerically that $\pi(2; x)-\pi(4; x)$ behaves as the uncorrelated random  walk;
also the number of returns of this random walk to the origin (the number of such
$x$ that  $\pi(2; x)=\pi(4; x)$) follows the usual square root  law
$\sqrt{x}$  \cite{Wolf1998}. It is in agreement with last sentences of the paper
\cite{Granville-Races}   suggesting that there is no Chebyshev bias for pairs of
primes ${p, p+2k}$ in general.

\vskip 0.4cm
\begin{center}
{\sf TABLE {\bf IV}}\\
\bigskip
\begin{tabular}{||c||c|c|c||c|c|c||} \hline
$x$ & $ \pi_q(x;3,2) $ & $ \pi_q(x;3,2)/\pi_q(x;3,1) $  & $ w_2(x)/\pi_q(x)$ & $\delta_1(x) $ & $ \delta_2(x) $ &  $\delta_0(x)  $ \\ \hline
$10^{3}$ &               8 &       4.000000  &        0.818182 &     0.000000 &     0.850210 &    0.000000 \\ \hline
$10^{4}$ &              12 &       1.714286  &        0.700000 &     0.023927 &     0.873320 &    0.024337 \\ \hline
$10^{5}$ &              32 &       1.684211  &        0.269231 &     0.217321 &     0.709371 &    0.030297 \\ \hline
$10^{6}$ &              71 &       1.731707  &        0.132743 &     0.340751 &     0.597523 &    0.028523 \\ \hline
$10^{7}$ &             215 &       2.128713  &        0.470032 &     0.341211 &     0.592423 &    0.039598 \\ \hline
$10^{8}$ &             560 &       1.992883  &        0.769596 &     0.301877 &     0.640202 &    0.034797 \\ \hline
$10^{9}$ &            1565 &       1.924969  &        0.322825 &     0.363780 &     0.582380 &    0.033401 \\ \hline
$10^{10}$ &            4410 &       1.963491  &        0.140454 &    0.424147 &     0.526904 &    0.030546 \\ \hline
$10^{11}$ &           12503 &       1.978636  &        0.049673 &    0.476509 &     0.479021 &    0.027770 \\ \hline
$10^{12}$ &           36069 &       1.999279  &        0.073571 &    0.514925 &     0.444198 &    0.025570 \\ \hline
$10^{13}$ &          104214 &       2.009254  &        0.671500 &    0.476702 &     0.485462 &    0.025570 \\ \hline
$10^{14}$ &          304482 &       2.004754  &        0.887649 &    0.442651 &     0.522215 &    0.023710 \\ \hline
$10^{15}$ &          893184 &       1.999557  &        0.610524 &    0.442717 &     0.524391 &    0.022016 \\ \hline
$10^{16}$ &         2636218 &       2.000222  &        0.709507 &    0.433234 &     0.535911 &    0.020650 \\ \hline
$10^{17}$ &         7817678 &       1.999806  &        0.474900 &    0.442119 &     0.528811 &    0.019378 \\ \hline
$10^{18}$ &        23265155 &       1.999574  &        0.169032 &    0.471831 &     0.500710 &    0.018267 \\ \hline
$10^{19}$ &        69497579 &       1.999852  &        0.141396 &    0.491985 &     0.481998 &    0.017256 \\ \hline
$10^{20}$ &       208240005 &       2.000040  &        0.365513 &    0.499533 &     0.475749 &    0.016352 \\ \hline
\end{tabular} \\
\end{center}
\vskip 0.4cm

All the primes of the form $q = 4k^2+1$ are necessarily disposed somewhere in the
arithmetical progression $4l+1$: $q\equiv 1 \pmod 4$.
But the situation changes when we consider residues modulo 3.
Let $\pi_q(x; 3, 1)$ denote the number of
primes $q=m^2+1<x$ such that $q \equiv 1 \pmod {3}$ and let  $\pi_q(x; 3, 2)$ denote
the number of primes $q=m^2+1<x$ such that $q \equiv 2 \pmod{3}$.  The direct
inspection of all possibilities shows that the residue 2  should appear twice
as often as the residue 1:
\bee
\pi_q(x; 3, 2)  \approx  2 \pi_q(x; 3, 1).
\eee
It translates into the human 10-base system as the observation, that except for
the first two cases 2 and 5,
the last digit of the primes of the form $m^2+1$ can be only 1 or 7,
see (\ref{set_Q}). The direct inspection of all 10 possibilities under
the assumption of normality of digits of $k$ in the base 10 shows,
that there are two times more ways  of obtaining 7 than 1.

The ratio of the number of those primes $p = 4k^2+1$ congruent to 2  modulo 3
to those congruent
to 1 modulo 3 can give some information about irregularities in the distribution
of primes of the form $4k^2+1$.
The Table IV shows the  values  of the ratio $\pi_q(x;3,2)/ \pi_q(x;3,1)$ for $x=10^3,\ldots 10^{20}$.
As it is seen from this Table after the initial transient interval below
 $10^7$ this ratio begins to
oscillate around the predicted value 2. Initially $\pi_q(x;3,2)>2 \pi_q(x;3,1)$
and for the first time $\pi_q(x;3,2) < 2 \pi_q(x;3,1)$ at $q_{17}=4\cdot 42^2+1 = 7057$,
i.e. at the 17-th prime of the form $q = m^2+1$. Recently A. Granville and G. Martin
\cite{Granville-Races}
have discussed several examples of "prime races". Primes of the form $q = 4k^2+1$
provide
another example of such a race. Namely we will say that two wins at a given $x$  if
$\pi_q(x;3,2) > 2 \pi_q(x;3,1)$ and let $w_2(x)$ denote the number of those primes
$q = 4k^2+1<x$ that residue two
wins over residue one. In the third column of the Table IV the ratio $w_2(x)/\pi_q(x)$
 is shown.
As it is  seen from this sample of numbers there are large fluctuations of the
ratio $w_2(x)/\pi_q(x)$. The Fig.7 shows the plot of $\pi_q(x;3,2)-2
 \pi_q(x;3,1)$
 up to $x=10^{20}$.  This plot can be interpreted as  a kind of one dimensional
random walk:   let $y(x)$ denote
the displacement of the walker at the ``time'' $x$, which plays the role of the time.
If for a given $q\in Q$ we find
that $ q   \equiv 1 \pmod 3 $ the random walker performs step down of length 2
and if  $ q \equiv  2  \pmod 3 $ the random walker performs step up of length 1
at the moment $x=q$. In other moments of time $x$ the walker simply  does not move.
Thus we have $y(x)= \pi_q(x;3,2) - 2 \pi_q(x;3,1)$.
This plot resembles usual random walk, there were 21349 returns to the origin of this
random walk up to $10^{20}$. There are large regions that $y(x)>0$ as well as
$y(x)<0$ suggesting that
there is no  Chebyshev bias in the distribution of primes $q_n$.
In the Fig.8 the plot of $y(x)/x^{1/4}=(\pi_q(x;3,2) - 2 \pi_q(x;3,1))/x^{1/4}$ is
shown. The amplitude of oscillations in this plot seems to be constant over the interval
$(10^3,  10^{20})$ and contained in the very short interval
$(-\frac{1}{2},~\frac{1}{2})$.
Dividing $(\pi_q(x;3,2) - 2 \pi_q(x;3,1))$  by $x^{\alpha}$
results in amplitudes going to zero when $\alpha>\frac{1}{4}$ and increasing when
$\alpha<\frac{1}{4}$. It is another argument in favor  of the error  term conjectured
in the Sect. 2.  There are probably logarithmic factors present,  like for the usual
Chebyshev bias, see Fig.  6 in \cite{Granville-Races}, but we are not able to
separate it. The amplitude of oscillations of $y(x)/x^{1/4}$ is very small,
less than 0.5 and roughly half of the plot in  Fig.8
is greater than zero and roughly half is below line zero. In \cite{Rubinstein_Sarnak} it was proposed
to  use the logarithmic density to  measure the Chebyshev bias. Here we will define
these densities for primes from $Q$ as follows:
\bee
\delta_1=\lim_{x \rightarrow \infty} \frac{1}{\log(x)}\sum_{\substack{2\leq n<x\\
2\pi_q(n;3,1)>\pi_q(n;3,2)}}\frac{1}{n}
\label{delta1}
\eee

\bee
\delta_2=\lim_{x \rightarrow \infty}\frac{1}{\log(x)}\sum_{\substack{2\leq n<x\\
2\pi_q(n;3,1)<\pi_q(n;3,2)}}\frac{1}{n}
\label{delta2}
\eee

\bee
\delta_0=\lim_{x \rightarrow \infty}\frac{1}{\log(x)}\sum_{\substack{2\leq n<x\\
2\pi_q(n;3,1)=\pi_q(n;3,2)}}\frac{1}{n}
\label{delta0}
\eee

We do not have at our disposal any  formulas like those in \cite{Rubinstein_Sarnak} and
we have to turn  to the brute force numerical calculation of finite size
approximations $\delta_1(x), \delta_2(x)$ and $\delta_0(x)$ given by expressions
(\ref{delta1})---
(\ref{delta0}) without limit operation $\lim_{x \rightarrow \infty}$.
The results are presented in
the 5-th, 6-th  and 7-th column of the Table IV
and in Figure 9. Up to $x=2^{31}$  the data for the Table IV and Figure 9 was obtained
by direct summing of the harmonic sums, for $x>2^{31}\approx 2.15\times 10^9$
the incredible accurate approximation   \cite{detemple1993}, \cite[pp. 76-78]{Havil03}
\bee
\sum_{k=n}^m \frac{1}{k}=\log\left(m+\frac{1}{2}\right)- \log\left(n-\frac{1}{2}\right) + {\mathcal{O}}\left(\frac{1}{n^2}\right)
\eee
was used, thus the error of each summand was smaller than $10^{-19}$, and as there
were $\mathcal{O}(10^8)$ terms, in the worst case of adding up all roundoffs  we
expect the total error to be smaller than $10^{-10}$. To calculate the harmonic
series up to $x=10^{20}$ directly by adding all numbers $1/n$ would take
a few thousands  years of CPU time, but it is in general impossible using standard
programming languages as the loop can have only integer counter and on 64 bit processors
the largest integer is $2^{63} \approx 9.22 \times 10^{18}$.

For all primes Chebyshev conjectured that
\bee
\lim_{x \rightarrow \infty} \sum_{p>2} (-1)^{\frac{p-1}{2}} e^{-p/x}=-\infty
\label{limit-F}
\eee
It was proved by Hardy and Littlewood \cite{H-L} and Landau \cite{Landau1}, \cite{Landau2}
that  (\ref{limit-F}) is true if and only if the L-function
\bee
L(s, \chi_1) = \sum_{n=0}^\infty \frac{(-1)^n}{(2n+1)^s}
\eee
does not have nontrivial zeros outside the critical line $\Re(s)=\frac{1}{2}$.
Here we formulate the  analogous function for the case of primes from $Q$:
\bee
F(x) = \sum_{q \in Q} c_q e^{-q/x}
\label{funkcja_F}
\eee
where
\bee
c_q\,=\,
\left\{
\begin{array}{ll}
2 & \mbox {when $q \!\!\!\! \mod 3 =1 $} \\
-1 & \mbox{\rm when  $q \!\!\!\! \mod 3 =2$}
\end{array}
\right.
\eee
We made the plot of $F(x)$ and in contrast to (\ref{limit-F}) this function seems to
not possess  a limit when $x \rightarrow \infty $. Because values of $F(x)$ sometimes
are close to zero calculating  the sum (\ref{funkcja_F}) on the computer we have
stopped summation at such $q'$  that
\bee
\left|\frac{c_{q'} e^{-q'/x}}{ \sum_{q=2}^{q'} c_q e^{-q/x}}\right|< 10^{-8},
\label{relative}
\eee
thus we have used relative error. Using such a condition (\ref{relative}) for terminating
the sum in (\ref{funkcja_F}) is necessary as $F(x)$ sometimes crosses zero and
using absolute error can be misleading at small values of $F(x)$.
The last 43  points in the Fig. 10 were
not fulfilling the requirement (\ref{relative}) as all primes  $q=m^2+1$  generated
were exhausted, but as it is seen in the plot of $F(x)$  in Fig. 10 values of the
function $F(x)$ are well above 1000 for $x>10^{17}$.

The numbers presented in the Table IV and plots in the Figures 8, 9 and 10 allow us to formulate the

{\bf Conjecture 5:} There is no Chebyshev bias for primes of the form $m^2+1$:
 $\delta_1=\delta_2=\frac{1}{2}, ~~~\delta_0=0$\\

It is the last conjecture formulated in this paper.

\bigskip


{\bf Acknowledgment:} I would like to thank Prof. Bernd Fischer for suggesting me in
December of 1996 the  search for primes of the form $m^2+1$. I thank Dr Jerzy
Cis{\l}o, Prof. Andrew Granville  and
Prof. W{\l}adys{\l}aw Narkiewicz for many comments and remarks.

\begin{figure}
\vspace{-0.3cm}
\begin{center}
\includegraphics[width=\textwidth,angle=0]{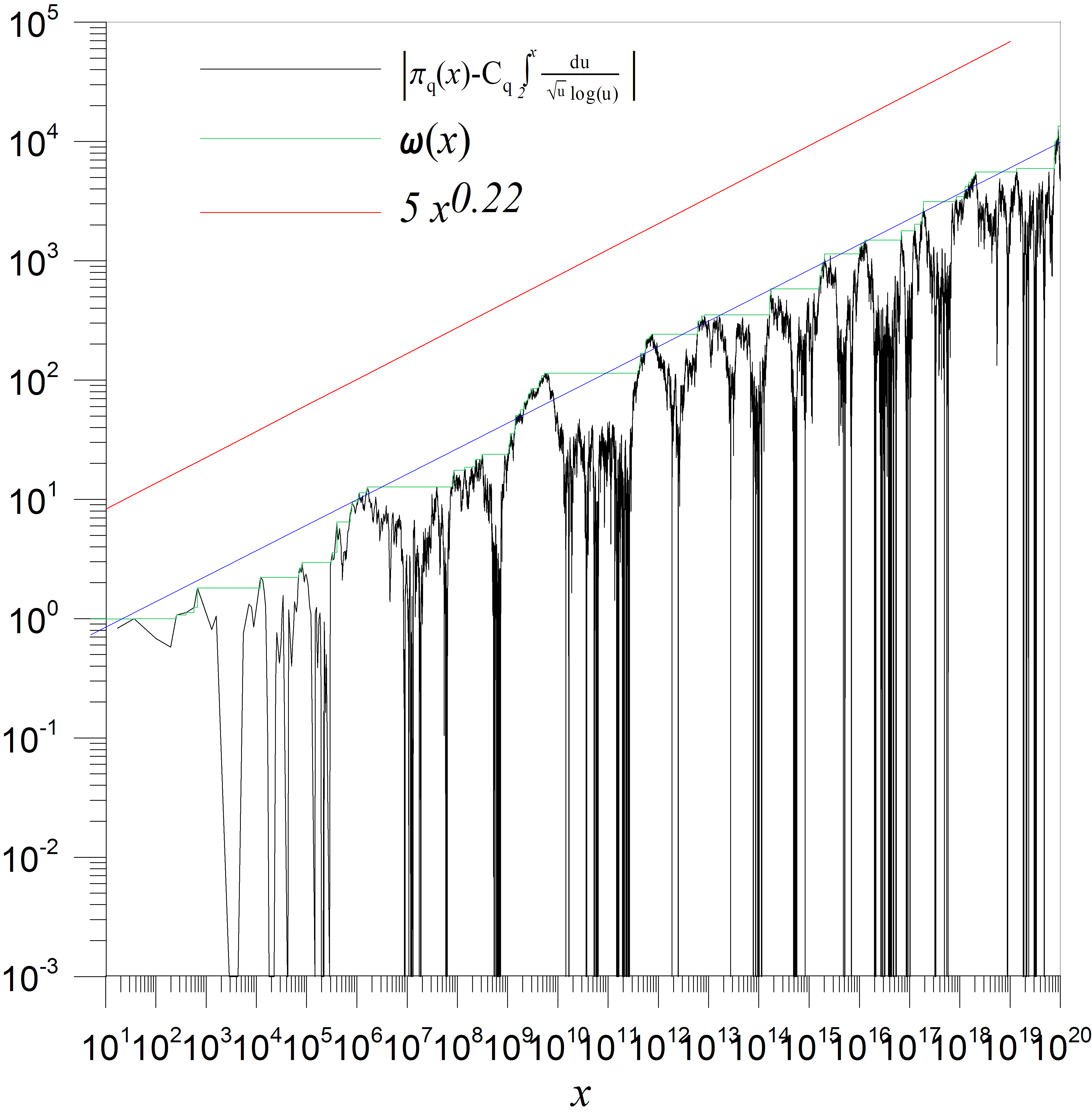} \\
\vspace{1.7cm} Fig.1 The plot of the error term for the conjecture (\ref{conj_calka})
up to $x=10^{20}$ on the double logarithmic axes. The plot of the difference
$\left| \pi_q(x) - {1 \over 2} C_q  \int_2^x {du \over \sqrt{u}\log{u}}\right|$ was
made from the 101566 points. Up to $10^{12}$  values of $|\Delta_q(x)|$ at all 54110 primes
of the form $m^2+1<10^{12}$ are included. For larger values of $x>10^{12}$
the following decimation procedure was used:  changes  of local maxima larger than 1,
all local minima and all sign changes were recorded with the provision that
values corresponding
to the sign changes of the $|\Delta_q(x)|$ and  smaller than $10^{-3}$ were artificially
set to $10^{-3}$. Additionally  the plot
contains 18430 points recorded at the geometrical progression $x=10^{12}\times (1.001)^n$.
In blue the  power-like fit $0.51877 \times x^{0.2139}$
 to $\omega(x)$ obtained by the least-square method  is
shown and in red the possible choice for the bound
of  $\omega(x)$ from above is plotted. \\
\end{center}
\end{figure}

\begin{figure}
\vspace{-0.9cm}
\begin{center}
\includegraphics[height=0.7\textheight,angle=90]{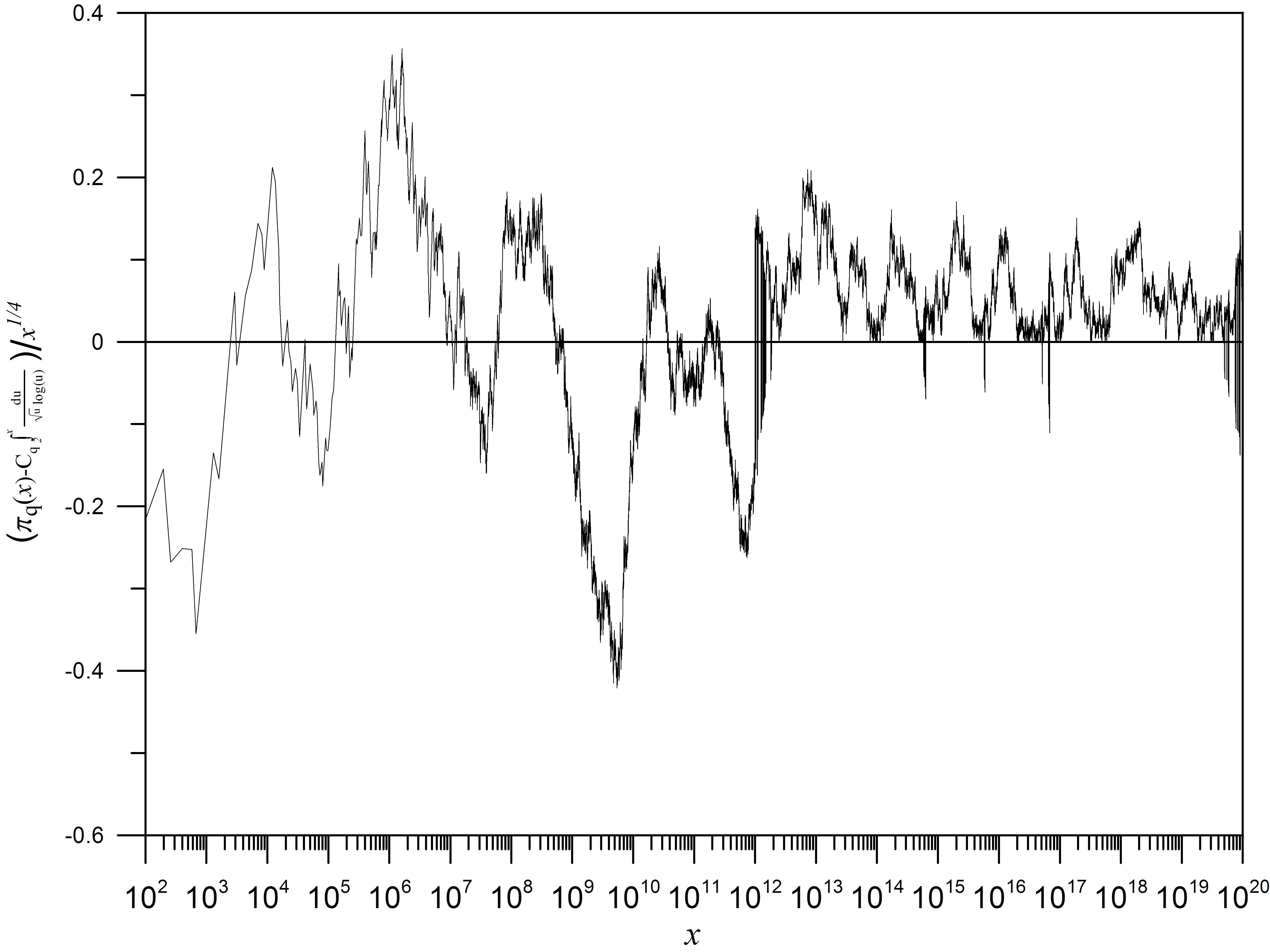} \\
\vspace{0.2cm} Fig.2 The plot of the oscillations of the
$(\pi_q(x) - {1 \over 2} C_q  \int_2^x {du \over \sqrt{u}\log{u}})/x^{1/4}$.
It consist of 106027 points and the procedure for recording values of this difference
was similar to the previous one.
\end{center}
\end{figure}

\begin{figure}
\vspace{-2.3cm}
\begin{center}
\includegraphics[width=\textwidth, angle=0]{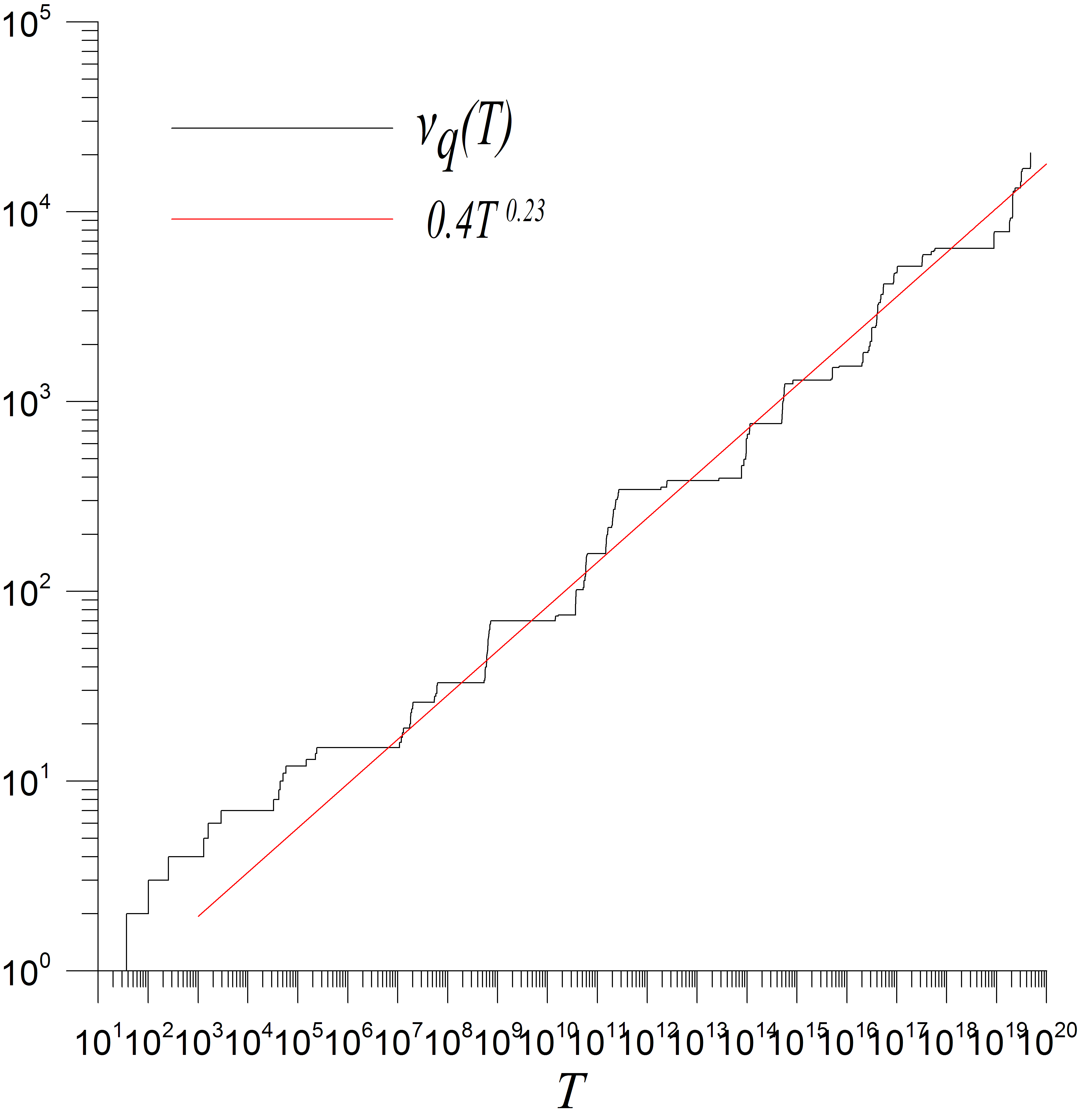} \\
\vspace{1.7cm} Fig.3 The plot of the number of sign changes of the difference
$\pi_q(x) - {1\over 2} C\int_2^x{du \over \sqrt{u}\log(u)}$ up to $x=10^{20}$. There are
20456 data points and power-like fit was  performed with respect to all points.
In red the power fit obtained by the least-square method is shown. \\
\end{center}
\end{figure}

\begin{figure}
\vspace{-1.0cm}
\begin{center}
\hspace{-1.5cm}
\includegraphics[height=0.5\textheight,angle=90]{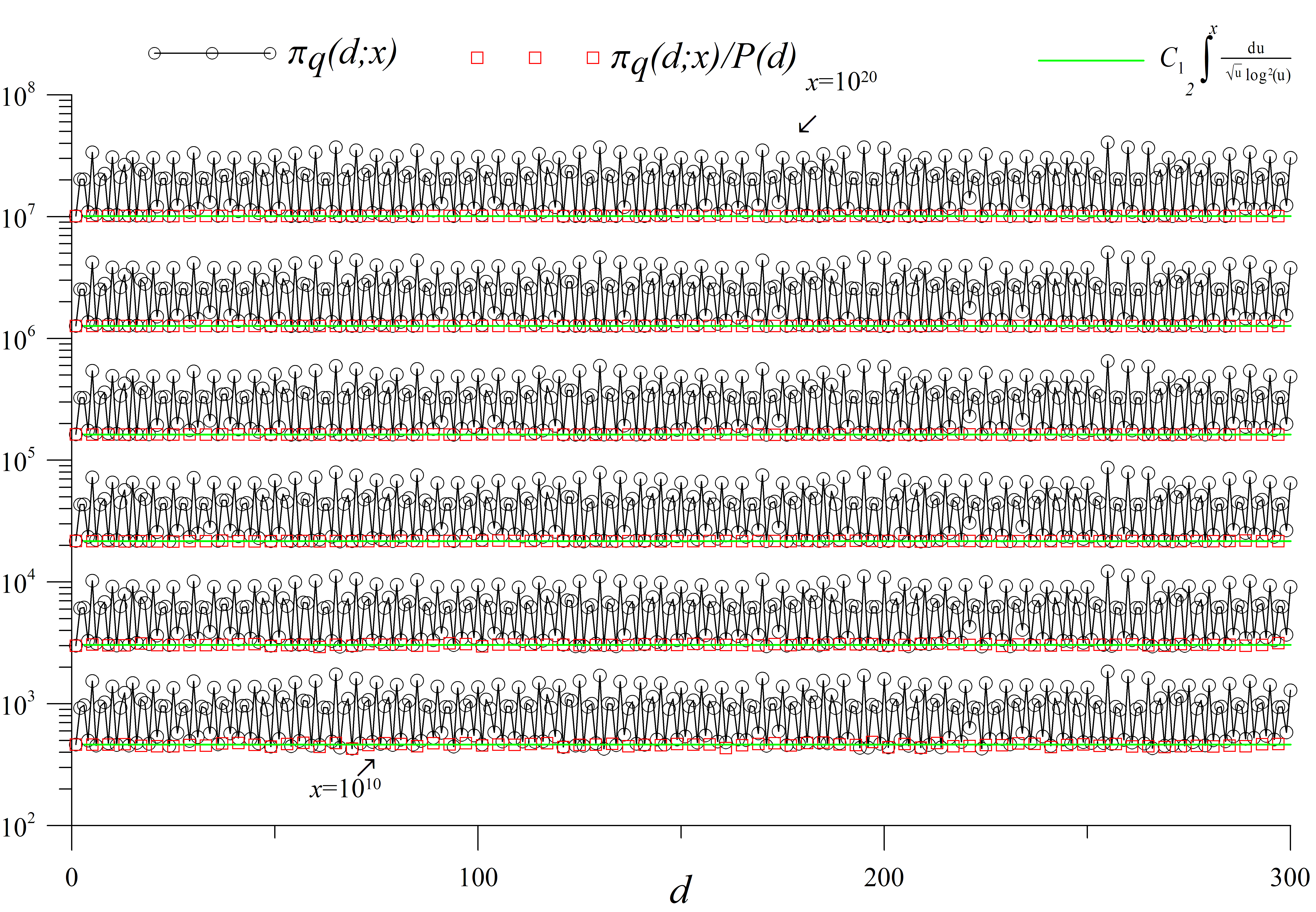} \\
\vspace{1.8cm} Fig.4 The plot of the function $\pi_q(d;x)$ for
$x=10^{10}, 10^{12}, \ldots, 10^{20}$ and $d\leq 300$. By red boxes the values of
$\pi_q(d;x)/P(d)$  are shown; for clarity only every 4-th value of $d$ is shown ---
descarding small fluctuations they fit the green lines representing the integrals
$C_1\int_5^x \frac{dt}{\sqrt{t}\log^2(t)}$.\\
\end{center}
\end{figure}

\begin{figure}
\vspace{-1.0cm}
\begin{center}
\hspace{-1.5cm}
\includegraphics[height=0.5\textheight,angle=90]{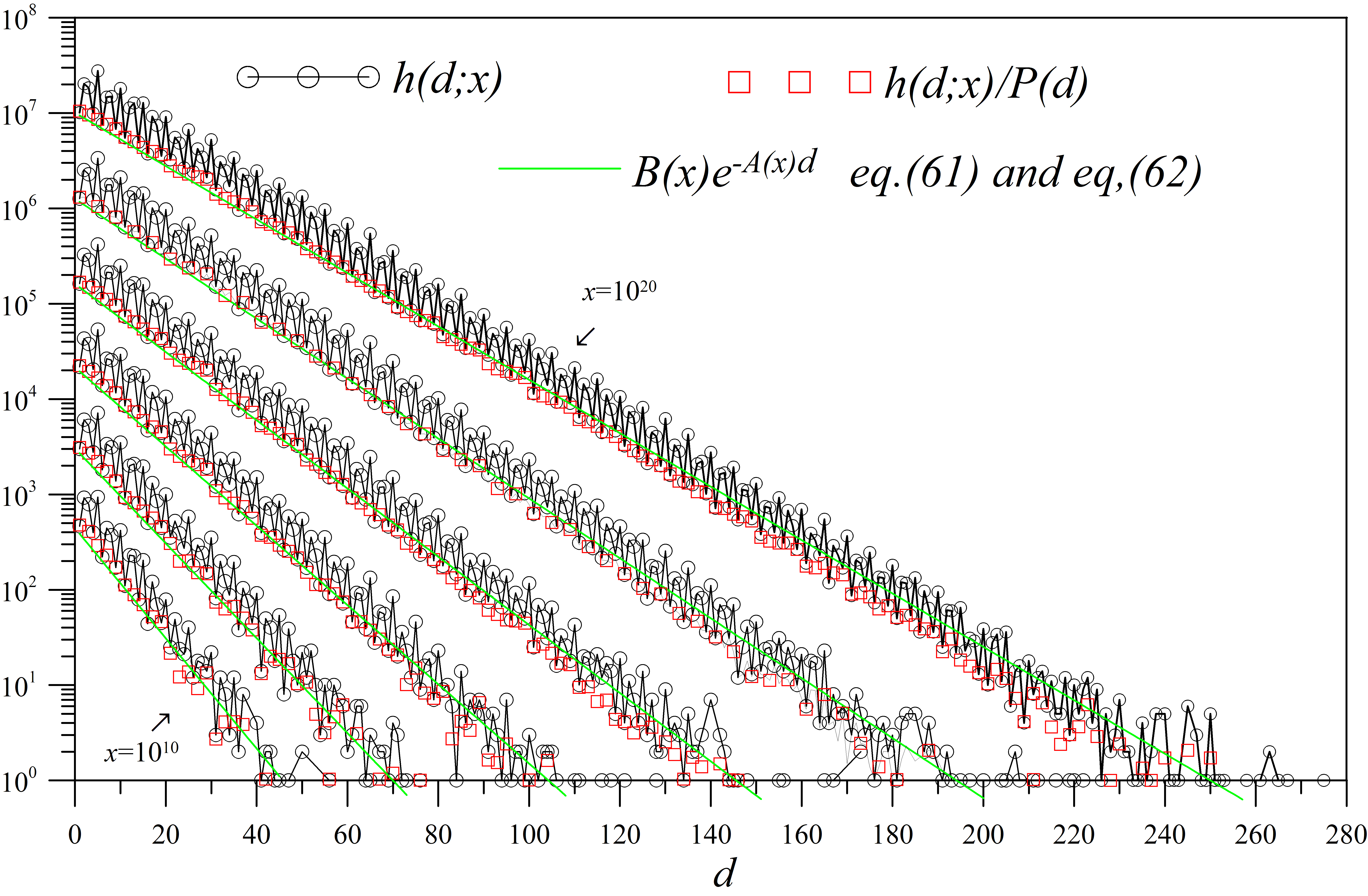} \\
\vspace{1.8cm} Fig.5 The plot of the histogram $h(d;x)$ for
$x=10^{10}, 10^{12}, \ldots, 10^{20}$. . By red boxes the values of
$h(x,d)/P(d)$ are shown; for clarity only every 4-th red box is shown. They
are supposed to lie along the green lines representing the plots of the
$\frac{\pi_q^2(x)}{s\sqrt{x}}\left(1-\frac{2\pi_q(x)}{\sqrt{x}}\right)^{d-1}$.
\end{center}
\end{figure}

\begin{figure}
\vspace{-2.3cm}
\hspace{-3.5cm}
\begin{center}
\includegraphics[width=\textwidth,angle=0]{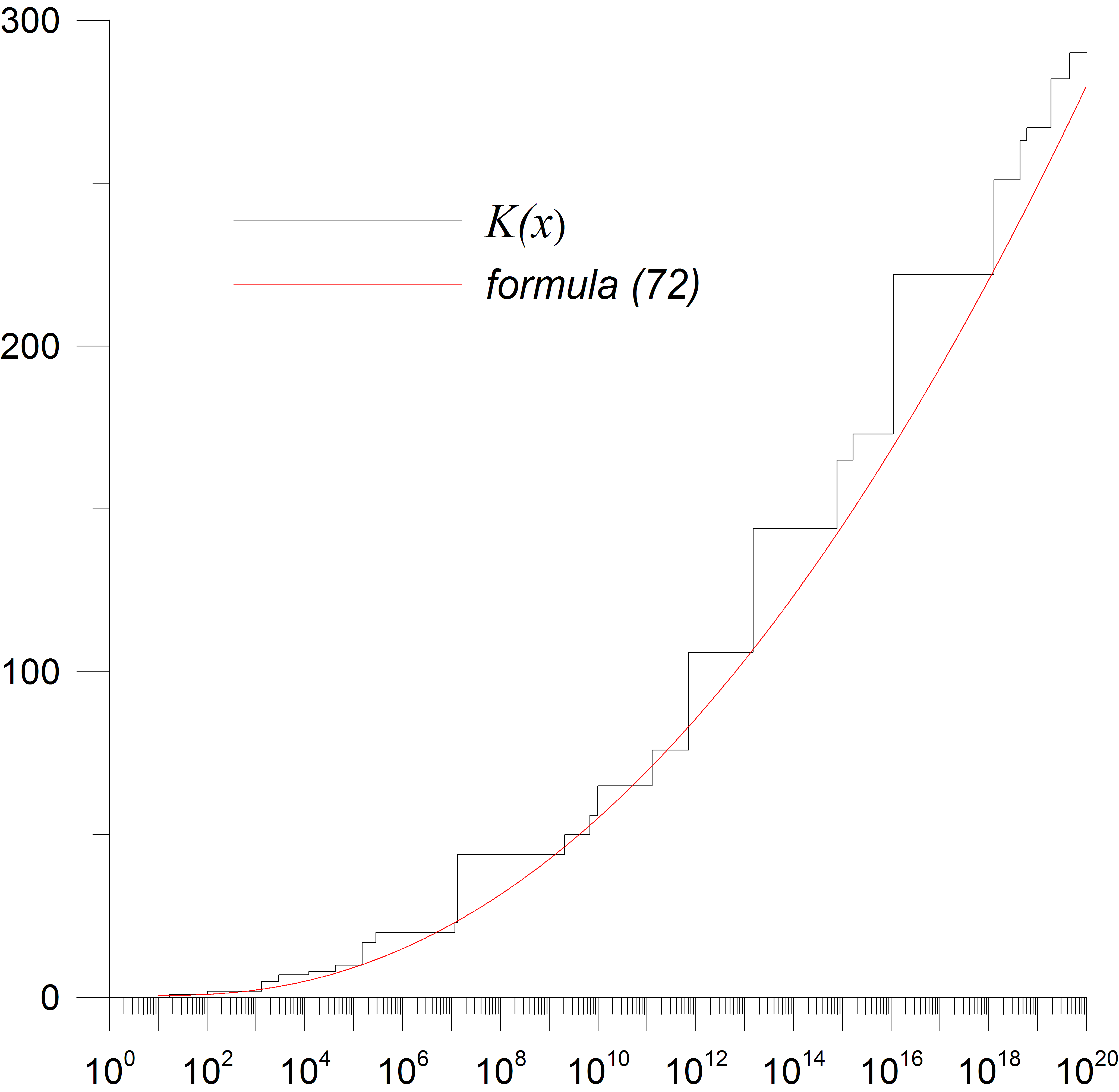} \\
\vspace{1.7cm} Fig.6 The plot of the function $K(x)$
\end{center}
\end{figure}

\begin{figure}
\vspace{-2.3cm}
\hspace{-3.5cm}
\begin{center}
\includegraphics[width=\textwidth,angle=0]{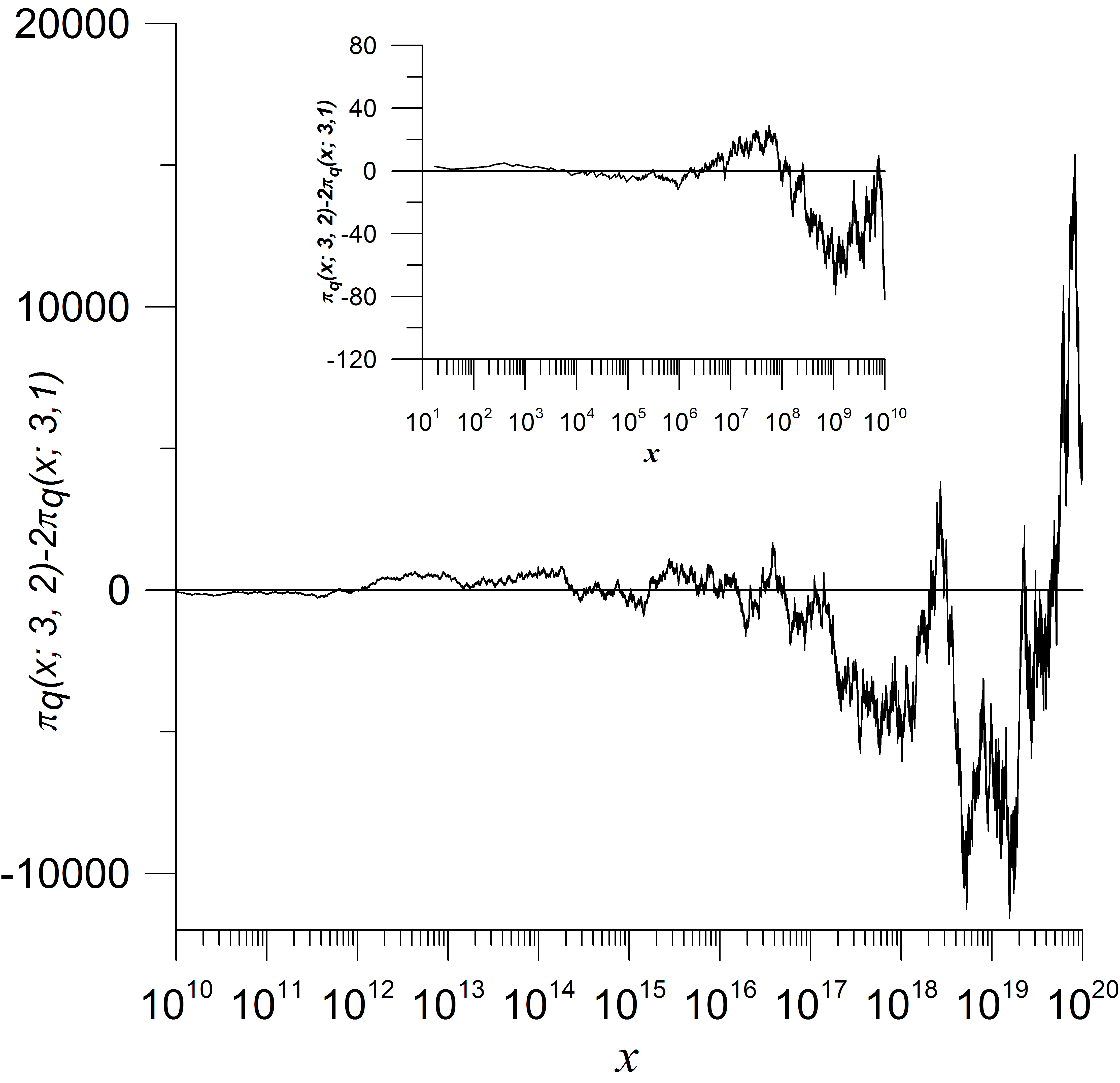} \\
\vspace{1.7cm} Fig.7 The plot of the difference $(\pi_q(x;3,2) - 2 \pi_q(x;3,1))$.
There are 102714 points: up to $10^{12}$ all primes $q_n$ are plotted, for $x>10^{12}$
the following decimations procedure was used: all local maxima, minima and sign changes
of this difference were recorded plus 18430 values at the geometrical progression
$10^{12}\times(1.001)^n$. Because the amplitude grows with $x$ we have presented
the interval $(10^{10}, 10^{20})$ in the inset  with smaller range on the $y$-axis.
\end{center}
\end{figure}

\begin{figure}
\vspace{-1.1cm}
\hspace{-3.5cm}
\begin{center}
\includegraphics[width=1.4\textwidth,angle=90]{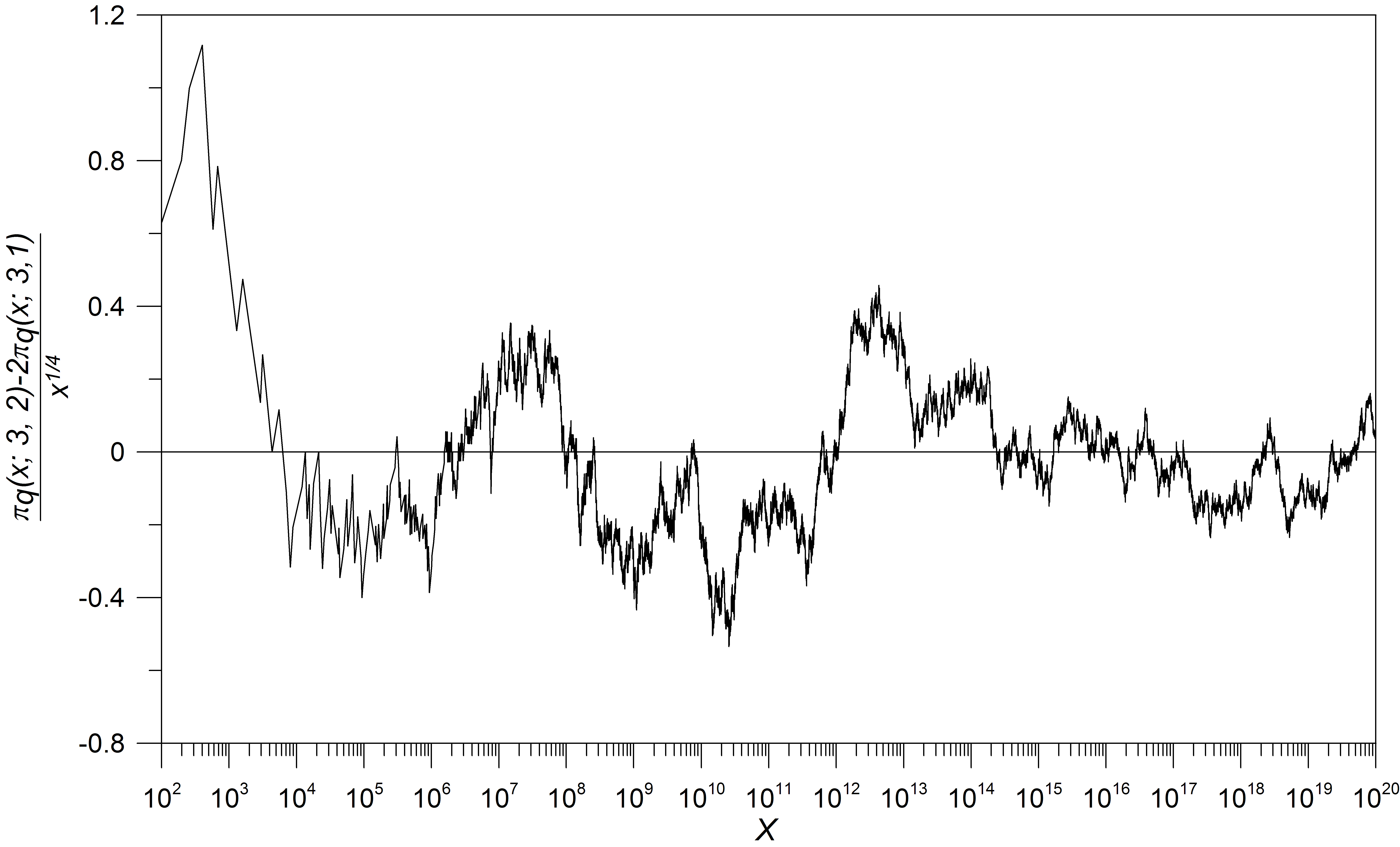} \\
\vspace{1.7cm} Fig.8 The plot of the ratio
 $y(x)/x^{1/4}=(\pi_q(x;3,2) - 2 \pi_q(x;3,1))/x^{1/4}$.
This plot is made of 102714 points, exactly as previous plot.
\end{center}
\end{figure}

\begin{figure}
\vspace{-2.3cm}
\hspace{-3.5cm}
\begin{center}
\includegraphics[width=\textwidth,angle=0]{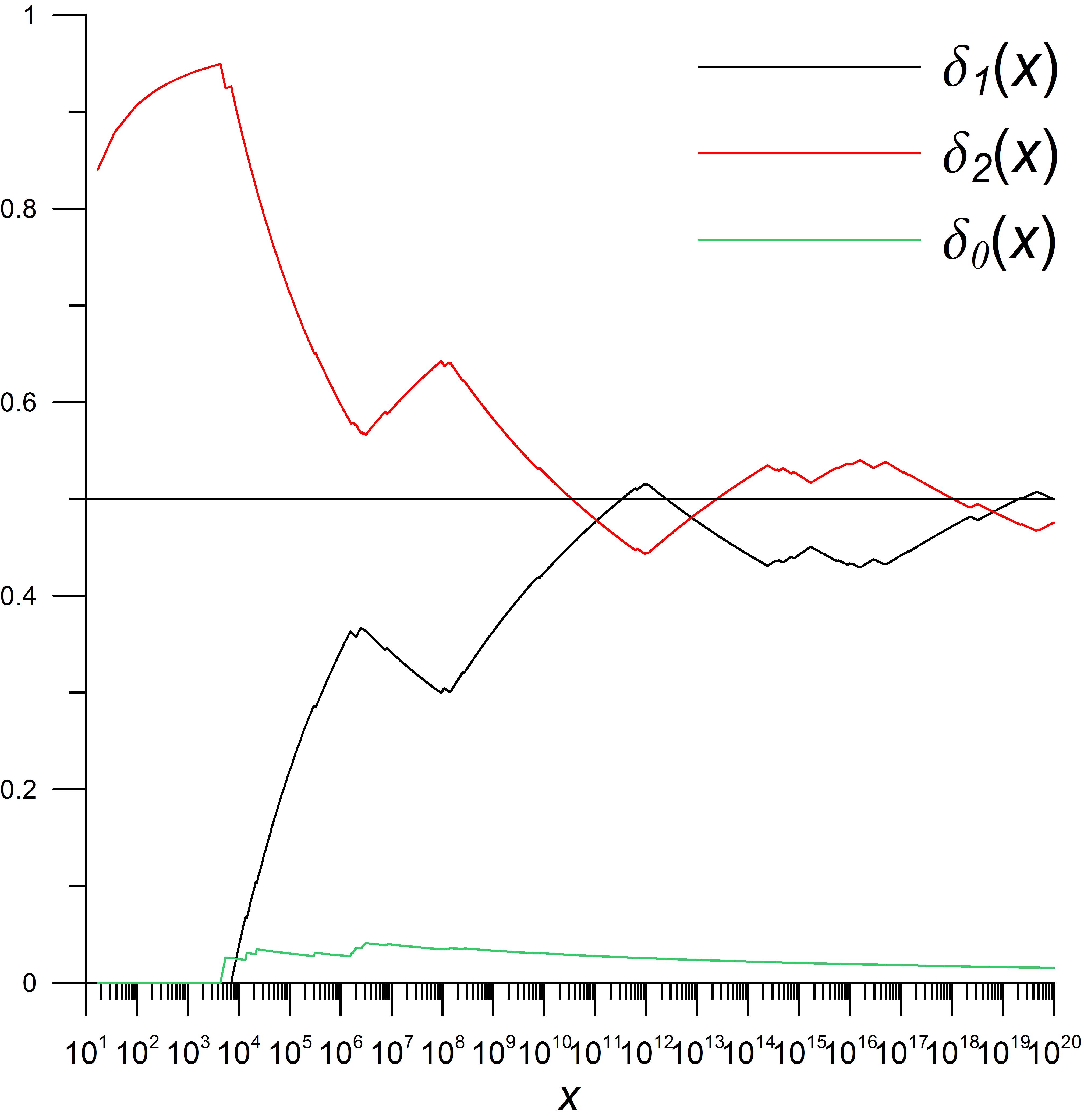} \\
\vspace{1.7cm} Fig.9 The plots of the logarithmic densities $\delta_1(x), \delta_2(x)$
and $\delta_0(x)$ defined in the text. Each plot consists of
72542 points: up to $10^{12}$  all  all primes $q_n$ are plotted, for $x>10^{12}$
the values of $\delta(x)$'s were  recorded at the progression
$10^{12}\times(1.001)^n$.
\end{center}
\end{figure}

\begin{figure}
\vspace{-1.0cm}
\begin{center}
\hspace{-1.5cm}
\includegraphics[height=0.5\textheight,angle=0]{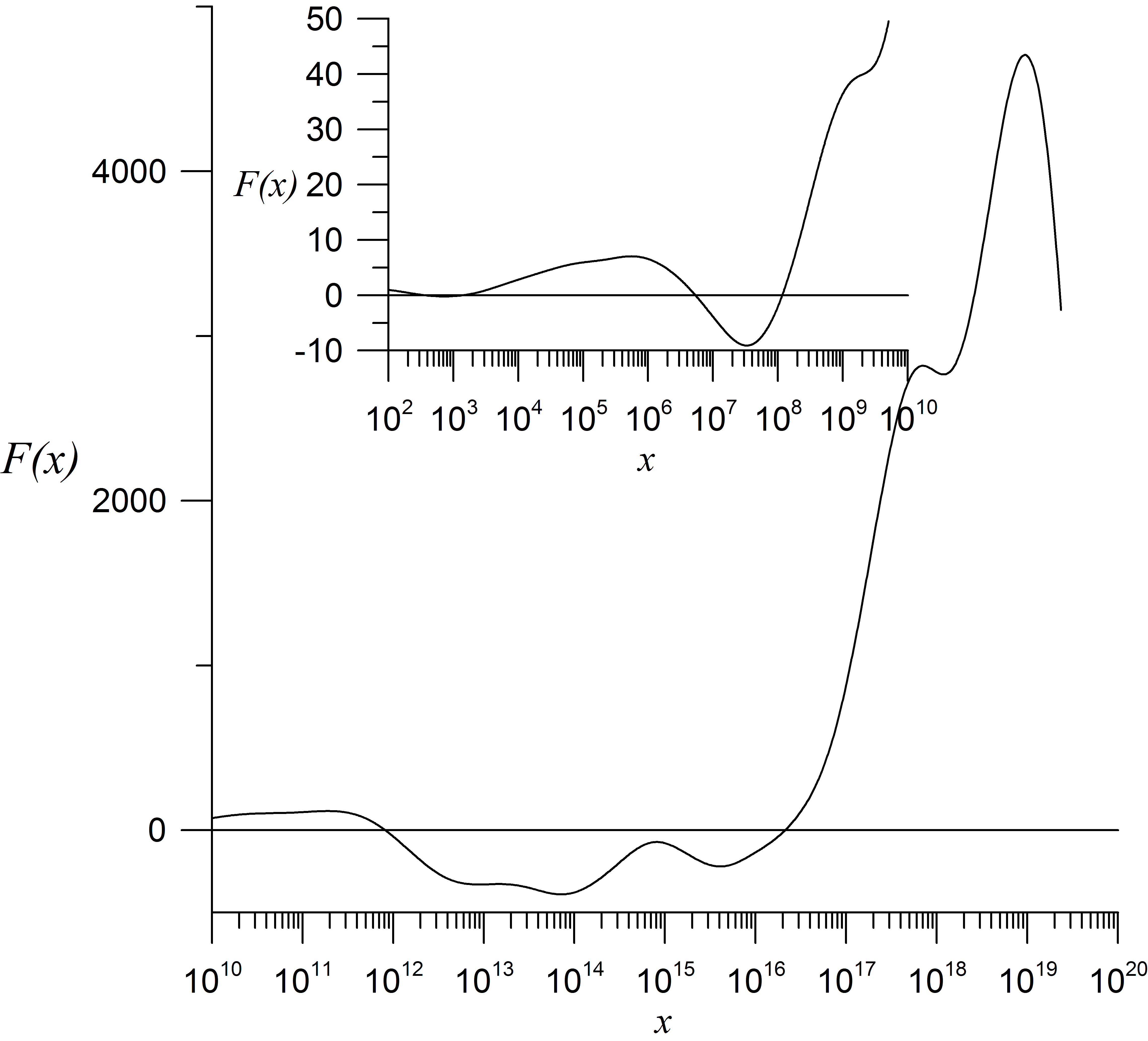} \\
\vspace{1.8cm} Fig.10 The plot of the function $F(x)$ defined by equation
(\ref{limit-F}). There are 2021 points on these figures: the arguments $x$ were
chosen as the geometrical progression $100.0\times 1.02 )^n$, $n=0, 1, \ldots, 2020$
and the largest value of $x$ is here  $ 2.35693149\ldots \times 10^{19}$, for which the
last term in the sum (\ref{funkcja_F}) was roughly $e^{-10^{20}/2.35693149\times 10^{19}}=0.01436723\ldots$.
\end{center}
\end{figure}

\end{document}